\documentclass[11pt]{amsart}

\usepackage{mathtools}
\usepackage{enumitem}
\usepackage{amsthm}
\usepackage{amsmath, amssymb}
\usepackage[utf8]{inputenc}
\usepackage[T1]{fontenc}
\usepackage{hyperref}

\theoremstyle{plain}
\newtheorem{Theorem}{Theorem}[section]
\newtheorem{Lemma}[Theorem]{Lemma}
\newtheorem{Proposition}[Theorem]{Proposition}
\newtheorem{Corollary}[Theorem]{Corollary}

\theoremstyle{definition}
\newtheorem{Definition}[Theorem]{Definition}
\newtheorem{Example}[Theorem]{Example}
\newtheorem{Remark}[Theorem]{Remark}

\hypersetup{
    pdftitle={From convolution pi-algebras to Rota-Baxter Hopf pi-algebras},
    pdfauthor={Xue Tang, Xin Wu, Yalan Zhang and Guodong Shi},
    pdfkeywords={Hopf pi-algebra; convolution pi-algebra; Rota-Baxter operator; antipode; group grading}
}

\title[Convolution $\pi$-algebras to Rota-Baxter Hopf $\pi$-algebras]{From convolution \(\pi\)-algebras to Rota-Baxter Hopf \(\pi\)-algebras}

\author{Xue Tang}
\address{Faculty of Mathematics and Physics, Huai'an University, Huai'an 223003, Jiangsu, China}

\author{Xin Wu}
\address{Faculty of Mathematics and Physics, Huai'an University, Huai'an 223003, Jiangsu, China}

\author{Yalan Zhang}
\address{Faculty of Mathematics and Physics, Huai'an University, Huai'an 223003, Jiangsu, China}

\author{Guodong Shi}
\address{Faculty of Mathematics and Physics, Huai'an University, Huai'an 223003, Jiangsu, China}
\email{guodongs@hau.edu.cn}

\date{} 

\begin{document}

\maketitle

\begin{abstract}
In this paper, we study Goncharov-type Rota--Baxter operators in the framework of group-graded Hopf \(\pi\)-algebras. By constructing convolution \(\pi\)-algebras as the core calculus tool, we systematically derive fundamental algebraic properties of Hopf \(\pi\)-algebras, including antipode anti-homomorphism and generalized homomorphism compatibility. On this basis, we establish the axiomatic system of two classes of graded Rota--Baxter operators of weight \(1\) and weight \(-1\), prove their bijective correspondence via the antipode, construct several non-trivial examples, and present two derived structures: descendant Hopf \(\pi\)-algebras and image Hopf \(\pi\)-subalgebras. This work enriches the operator theory of graded Hopf algebras and provides algebraic constructions for related research on graded Yang--Baxter equations.
\end{abstract}


\section{Introduction}

Rota-Baxter operators and Hopf algebras are core structures in algebra with profound theoretical implications and broad applications. The Rota-Baxter operator was first introduced by Baxter in 1960 to characterize the Spitzer identity in stochastic analysis. Later, it was systematically axiomatized by Rota, and Cartier investigated the structure of free Rota-Baxter algebras, clarifying the intrinsic connections between such operators, shuffle products and multiple zeta values. Since then, it has developed into an independent algebraic branch \cite{baxter1960analytic,rota1969baxter,cartier1972structure}. In mathematical physics, Connes and Kreimer revealed the deep relation between Hopf algebras and the renormalization procedure in quantum field theory, which provided the core motivation for the interdisciplinary study of these two structures \cite{connes2000renormalization}. Since the beginning of the 21st century, Guo and his collaborators have systematically combined Rota-Baxter operators with classical Hopf algebras, greatly advancing the development of the field. His monograph gives a comprehensive summary of the fundamental theories in this direction \cite{andrews2003baxter,guo2012introduction}.

Hopf algebras originated from Hopf’s topological research on Lie group manifolds. The algebraic foundations were gradually laid by Borel, Cartier, Milnor and Moore, and the modern axiomatic system together with Sweedler’s notation was finally established by Sweedler. He introduced core tools including convolution products and convolution inverses, making Hopf algebras an independent algebraic discipline separated from topology \cite{hopf1941topologie,milnor1965structure,sweedler1969hopf}. The theoretical framework of modern Hopf algebras was further refined in Radford’s textbook \cite{radford2012hopf}. In 2000, Turaev introduced Hopf \(\pi\)-algebras and Hopf \(\pi\)-coalgebras (collectively called Hopf group (co)algebras) indexed by groups, extending classical Hopf algebras to the group-graded setting. Hopf \(\pi\)-coalgebras can be used to construct invariants of principal \(\pi\)-bundles and 3-manifolds, and their representation categories naturally form crossed \(\pi\)-categories. When the indexing group is commutative, such structures reduce to \(\pi\)-colored ribbon Hopf algebras in the sense of Ohtsuki \cite{turaev2000homotopy,turaev2010homotopy,ohtsuki1993colored}. Virelizier further systematically established the basic theories of Hopf \(\pi\)-coalgebras, including integrals, traces and quasitriangular structures, laying a rigorous algebraic foundation for their topological applications \cite{virelizier2002hopf}. It should be noted that existing foundational research on group-graded Hopf structures mainly focuses on the coalgebra side, and the core properties related to convolution operations and antipode operators in the framework of Hopf \(\pi\)-algebras have not been systematically investigated.

In 2021, Goncharov established a complete theory of Rota-Baxter operators on cocommutative Hopf algebras, constructed the correspondence among group structures, Lie algebras and Rota-Baxter operators, and revealed the intrinsic relation between these operators and solutions of the Yang–Baxter equation \cite{goncharov2021rota}. In the group-graded direction, Ma et al. introduced classical weight-\(\lambda\) Rota-Baxter operators into the framework of Hopf group (co)algebras in 2022, and established the theory of Rota-Baxter Hopf T-algebras, covering graded versions of Rota-Baxter family algebras and Rota-Baxter bialgebras \cite{ma2022rotabaxter}. This work follows the classical definitional paradigm of Rota-Baxter operators on associative algebras, and belongs to the graded extension of the classical Rota-Baxter system. It belongs to a different research line from Goncharov’s Hopf-type Rota-Baxter operators, and the two definitional systems are incompatible with each other.

The aim of this paper is to systematically extend the Rota-Baxter operator theory on cocommutative Hopf algebras established by Goncharov to the framework of Hopf \(\pi\)-algebras, and to establish a group-graded version that is fully compatible with the classical Hopf-type results. During the derivation, we find that this extension cannot be directly based on existing literature: due to the lack of systematic conclusions on basic properties such as convolution operations and antipode compatibility on the Hopf \(\pi\)-algebra side, we need to build the required theoretical framework step by step starting from basic tools.

Given that group-graded Hopf structures involve multiple group indices, families of graded multiplications, graded comultiplications, graded antipodes and multi-layer composite maps, the notational system is hierarchically complex and prone to index confusion and mapping ambiguity. Section 2 is devoted to establishing a unified set of notational conventions throughout the paper, covering basic isomorphisms of vector spaces, layered naming rules for \(\pi\)-algebras, coalgebras and Hopf \(\pi\)-algebras, as well as the graded version of the Heyneman–Sweedler notation. This notational system is fully self-consistent within the scope of this paper. Although it contains a relatively large number of entries, it fundamentally standardizes the use of notations in subsequent derivations and ensures the rigor and readability of the proofs. Section 3 collects the basic definitions and preliminaries of Hopf \(\pi\)-algebras, including the rigorous axioms of \(\pi\)-algebras, Hopf \(\pi\)-algebras, generalized graded homomorphisms, as well as basic constructions such as opposite \(\pi\)-algebras, serving as the starting point of the whole paper. Section 4 constructs convolution \(\pi\)-algebras induced by coalgebras and \(\pi\)-algebras, verifies their graded algebra axioms, and proves that coalgebra homomorphisms and \(\pi\)-algebra homomorphisms can respectively induce homomorphisms of convolution \(\pi\)-algebras, so as to establish a complete graded convolution calculus. Using this as the core tool, we derive convolution inverse formulas for two typical classes of maps, and further obtain the antipode compatibility criterion for generalized homomorphisms of Hopf \(\pi\)-algebras. We also prove that the antipode operator is an anti-homomorphism of graded algebras and an anti-homomorphism of coalgebras, and establish the bijectivity of the antipode and the Hopf \(\pi\)-algebra structure of the coopposite graded structure in the cocommutative case. These properties are necessary prerequisites for constructing graded Rota-Baxter operators, and have not been systematically discussed in existing literature on Hopf \(\pi\)-algebras. On this basis, Section 5 establishes the axiomatic system of Rota-Baxter operators on cocommutative Hopf \(\pi\)-algebras that is fully compatible with Goncharov’s system, gives the definitions of two families of graded operators of weight 1 and weight \(-1\), and proves that they form a bijective correspondence via the antipode. We construct several explicit examples to verify the existence and richness of the structure, and discuss the basic properties of the operators as well as related derived structures.

\section{Notations}

In this section, we establish the notational conventions adopted throughout the paper. The primary purpose of these uniform conventions is to simplify the exposition and avoid any notational ambiguity. A precise understanding of the symbols defined below is essential for following the arguments and results presented in subsequent sections.
    
Throughout this paper, we work over a ground field \(\mathbb{K}\) of characteristic zero, unless otherwise specified. The symbol \(\pi\) always denotes an arbitrary group serving as the index set, and \(p, q, r\) always stand for arbitrary elements of \(\pi\).

For any vector space \(V\) over \(\mathbb{K}\), we fix the following notation:
\begin{enumerate}[label=(\roman*)]
\item \(\operatorname{id}_V\) denotes the identity map on \(V\);
\item \(V^*\) denotes its dual space \(\operatorname{Hom}_{\mathbb{K}}(V,\mathbb{K})\);
\item \(\phi\) together with its inverse \(\phi^{-1}\) denote the canonical linear isomorphism between \(V\) and \(\mathbb{K}\otimes V = \{k\otimes v \mid k\in\mathbb{K},\, v\in V\}\), given by \(\phi(v) = 1_{\mathbb{K}}\otimes v\) and \(\phi^{-1}(k\otimes v) = k\cdot v\) for all \(v\in V\) and \(k\in\mathbb{K}\);
\item \(\psi_V\) together with its inverse \(\psi_V^{-1}\) denote the canonical linear isomorphism between \(\operatorname{Hom}_{\mathbb{K}}(\mathbb{K}, V)\) and \(V\), given by \(\psi_V(f) = f(1_{\mathbb{K}})\) and \(\big(\psi_V^{-1}(v)\big)(k) = k \cdot v\) for all \(f \in \operatorname{Hom}_{\mathbb{K}}(\mathbb{K}, V)\), \(v \in V\) and \(k \in \mathbb{K}\).
\end{enumerate}

The notation conventions for families of graded algebraic structures below are used consistently throughout the paper.
\begin{enumerate}[label=(\alph*)]
\item A \(\pi\)-algebra is denoted by
\[
H_{\pi}:=\bigl(\{H_{p}\}_{p \in \pi},\{m_{p, q}\}_{p, q \in \pi}, \mu\bigr).
\]

We use superscripts to distinguish distinct \(\pi\)-algebras sharing the same index set. For example, two distinct \(\pi\)-algebras are given by
\[
H_{\pi}^{1}:=\bigl(\{H_{p}^{1}\}_{p \in \pi},\{m_{p, q}^{1}\}_{p, q \in \pi}, \mu^{1}\bigr),\quad H_{\pi}^{2}:=\bigl(\{H_{p}^{2}\}_{p \in \pi},\{m_{p, q}^{2}\}_{p, q \in \pi}, \mu^{2}\bigr).
\]
When the index set is replaced by \(\pi_{1}\) or \(\pi_{2}\), the corresponding algebras are denoted respectively by
\[
H_{\pi_{1}}:=\bigl(\{H_{p}\}_{p \in \pi_{1}},\{m_{p, q}\}_{p, q \in \pi_{1}}, \mu\bigr),\quad H_{\pi_{2}}:=\bigl(\{H_{p}\}_{p \in \pi_{2}},\{m_{p, q}\}_{p, q \in \pi_{2}}, \mu\bigr).
\]
Superscripts for distinguishing distinct algebras are independent of the choice of index set. For example, two distinct \(\pi_{1}\)-algebras are written as
\[
H_{\pi_{1}}^{1}:=\bigl(\{H_{p}^{1}\}_{p \in \pi_{1}},\{m_{p, q}^{1}\}_{p, q \in \pi_{1}}, \mu^{1}\bigr),\quad H_{\pi_{1}}^{2}:=\bigl(\{H_{p}^{2}\}_{p \in \pi_{1}},\{m_{p, q}^{2}\}_{p, q \in \pi_{1}}, \mu^{2}\bigr).
\]

We also use different capital base letters to refer to different families of \(\pi\)-algebras. When multiple base letters appear in the same context, the multiplication morphisms and unit morphism are labelled with the corresponding base letter to avoid ambiguity. For example, a \(\pi\)-algebra with base letter \(A\) takes the form
\[
A_{\pi}:=\bigl(\{A_{p}\}_{p \in \pi},\{m_{p, q}^{A}\}_{p, q \in \pi}, \mu^{A}\bigr).
\]
All preceding notation rules for superscripts and index sets apply identically to \(A_{\pi}\) and any other letter-labelled \(\pi\)-algebras.

\item A coassociative counital coalgebra is denoted by 
\[
\mathcal{C} := (C, \Delta, \varepsilon).
\] 
We use different calligraphic capital letters to denote distinct coalgebras; for example, a second coalgebra is written as \[
\mathcal{D} := (D, \Delta_{D}, \varepsilon_{D}).
\]
 We also use subscript indices to refer to indexed families of coalgebras. For instance, a family of coalgebras indexed by \(\pi\) takes the form
\[
\{\mathcal{H}_{p}\}_{p \in \pi} = \{(H_{p}, \Delta_{p}, \varepsilon_{p})\}_{p \in \pi}.
\]
When there are multiple such families indexed by \(\pi\), they are distinguished by superscripts. For example, two distinct families of coalgebras indexed by \(\pi\) are given by
\[
\{\mathcal{H}_{p}^{1}\}_{p \in \pi} := \{(H_{p}^{1}, \Delta_{p}^{1}, \varepsilon_{p}^{1})\}_{p \in \pi},\quad \{\mathcal{H}_{p}^{2}\}_{p \in \pi} := \{(H_{p}^{2}, \Delta_{p}^{2}, \varepsilon_{p}^{2})\}_{p \in \pi}.
\]
We adopt the Heyneman--Sweedler notation for comultiplications throughout this paper. 
For a coalgebra \((C, \Delta, \varepsilon)\) and any element \(c \in C\), we write
$
\Delta(c) = c_{(1)} \otimes c_{(2)}.
$
For the indexed family of coalgebras \(\{(H_p, \Delta_p, \varepsilon_p)\}_{p\in\pi}\), given any \(p \in \pi\) and any \(h \in H_p\), we write
$
\Delta_p(h) = h_{(1,p)} \otimes h_{(2,p)}.$

\item A Hopf \(\pi\)-algebra is denoted by
\[
\mathcal{H}_{\pi}:=\bigl(\{H_{p}, \Delta_{p}, \varepsilon_{p}\}_{p \in \pi},\{m_{p, q}\}_{p, q \in \pi}, \mu, \{S_{p}\}_{p \in \pi}\bigr).
\]
We use superscripts to distinguish distinct Hopf \(\pi\)-algebras sharing the same index set. For example, two distinct Hopf \(\pi\)-algebras are given by
\[
\mathcal{H}_{\pi}^{1}:=\bigl(\{H_{p}^{1}, \Delta_{p}^{1}, \varepsilon_{p}^{1}\}_{p \in \pi},\{m_{p, q}^{1}\}_{p, q \in \pi}, \mu^{1}, \{S_{p}^{1}\}_{p \in \pi}\bigr),
\]
\[
\mathcal{H}_{\pi}^{2}:=\bigl(\{H_{p}^{2}, \Delta_{p}^{2}, \varepsilon_{p}^{2}\}_{p \in \pi},\{m_{p, q}^{2}\}_{p, q \in \pi}, \mu^{2}, \{S_{p}^{2}\}_{p \in \pi}\bigr).
\]

When the index set is replaced by \(\pi_{1}\) or \(\pi_{2}\), the corresponding Hopf algebras are denoted respectively by
\[
\mathcal{H}_{\pi_{1}}:=\bigl(\{H_{p}, \Delta_{p}, \varepsilon_{p}\}_{p \in \pi_{1}},\{m_{p, q}\}_{p, q \in \pi_{1}}, \mu, \{S_{p}\}_{p \in \pi_{1}}\bigr),
\]
\[
\mathcal{H}_{\pi_{2}}:=\bigl(\{H_{p}, \Delta_{p}, \varepsilon_{p}\}_{p \in \pi_{2}},\{m_{p, q}\}_{p, q \in \pi_{2}}, \mu, \{S_{p}\}_{p \in \pi_{2}}\bigr).
\]
Superscripts for distinguishing distinct algebras are independent of the choice of index set. For example, two distinct Hopf \(\pi_{1}\)-algebras are written as
\[
\mathcal{H}_{\pi_{1}}^{1}:=\bigl(\{H_{p}^{1}, \Delta_{p}^{1}, \varepsilon_{p}^{1}\}_{p \in \pi_{1}},\{m_{p, q}^{1}\}_{p, q \in \pi_{1}}, \mu^{1}, \{S_{p}^{1}\}_{p \in \pi_{1}}\bigr),
\]
\[
\mathcal{H}_{\pi_{1}}^{2}:=\bigl(\{H_{p}^{2}, \Delta_{p}^{2}, \varepsilon_{p}^{2}\}_{p \in \pi_{1}},\{m_{p, q}^{2}\}_{p, q \in \pi_{1}}, \mu^{2}, \{S_{p}^{2}\}_{p \in \pi_{1}}\bigr).
\]

\item A Rota-Baxter Hopf \(\pi-\)algebra is denoted by 
\[
\prescript{R}{}{\mathcal{H}}_{\pi}:=(\{H_p, \Delta_p, \varepsilon_p\}_{p \in \pi}, \{m_{p,q}\}_{p,q \in \pi}, \mu, \{S_p\}_{p \in \pi}, \{B_p\}_{p \in \pi})
\]
\end{enumerate}

\section{Preliminaries}
	
This section collects the preliminary knowledge required throughout the paper, including the axiomatic definitions of graded algebraic structures and several basic lemmas. These results serve as a necessary foundation for the subsequent construction of convolution \(\pi\)-algebras and the study of Rota--Baxter operators.

We first recall the rigorous definitions of \(\pi\)-algebras, Hopf \(\pi\)-algebras and related basic concepts.

\begin{Definition}
\begin{enumerate}[label=(\alph*)]
\item A graded algebra indexed by \(\pi\) (abbreviated as \(\pi\)-algebra) \(H_\pi \) consists of a family of vector spaces \(\{H_p\}_{p \in \pi}\), 
a family of bilinear operator \(\{m_{p,q}: H_p \otimes H_q \to H_{pq}\}_{p,q \in \pi}\)
and a linear operator \(\mu: \mathbb{K} \to H_e\) such that for all \(p, q, r \in \pi\),
\begin{equation*}
\begin{gathered}
m_{pq,r}\circ(m_{p,q} \otimes \mathrm{id}_{H_r}) = m_{p,qr}\circ(\mathrm{id}_{H_p} \otimes m_{q,r}), \\
m_{p,e}\circ(\mathrm{id}_{H_p} \otimes \mu)\circ\phi = \mathrm{id}_{H_p} = m_{e,p}\circ(\mu \otimes \mathrm{id}_{H_p})\circ\phi.
\end{gathered}
\end{equation*}

\item A graded Hopf algebra indexed by \(\pi\) (abbreviated as Hopf \(\pi\)-algebra) 
\(\mathcal{H}_\pi\) 
consists of a \(\pi\)-algebra 
 \(H_\pi \) 
together with a family of coassociative counital coalgebras 
\(\{\mathcal{H}_p\}_{p \in \pi}\)
and a family of linear operators 
\(\{S_p : H_p \to H_{p^{-1}}\}_{p \in \pi}\)such that each multiplication map \(m_{p,q}\) and the unit \(\mu\) are coalgebra homomorphisms, i.e.
\begin{equation*}
\begin{gathered}
\Delta_{pq} \circ m_{p,q} = \bigl( m_{p,q} \otimes m_{p,q} \bigr) \circ \sigma_{23} \circ \bigl( \Delta_p \otimes \Delta_q \bigr), \quad m_{\mathbb{K}} \circ (\varepsilon_p \otimes \varepsilon_q) = \varepsilon_{pq} \circ m_{p,q}, \\
\Delta_e \circ \mu = (\mu \otimes \mu) \circ \Delta_{\mathbb{K}}, \quad \varepsilon_e \circ \mu = \mathrm{id}_{\mathbb{K}}
\end{gathered}
\end{equation*}
and for every \(p \in \pi\),
\begin{equation*}
m_{p^{-1},p}\circ(S_p \otimes \mathrm{id}_{H_p})\circ\Delta_p
= \mu \circ \varepsilon_p
= m_{p,p^{-1}}\circ(\mathrm{id}_{H_p} \otimes S_p)\circ\Delta_p.
\end{equation*}

\item Let \(\varPhi \colon \pi_1 \to \pi_2\) be a group homomorphism (isomorphism).
\begin{enumerate}[label=(\roman*)]
\item A \emph{generalized graded algebra homomorphism (isomorphism)} from a \(\pi_1\)-algebra \(H_{\pi_1}^1\) to a \(\pi_2\)-algebra \(H_{\pi_2}^2\) is a family of linear maps \(\{ f_p \colon H_p^1 \to H_{\varPhi(p)}^2 \}_{p\in\pi_1}\) such that for all \(p, q \in \pi_1\),
    \[
    m_{\varPhi(p),\varPhi(q)}^2 \circ(f_p \otimes f_q) = f_{pq}\circ m_{p,q}^1, \quad f_e\circ\mu^1 = \mu^2.
    \]
\item A \emph{generalized graded  Hopf algebra homomorphism (isomorphism)} from a Hopf \(\pi_1\)-algebra \(\mathcal{H}_{\pi_1}^1\) to a Hopf \(\pi_2\)-algebra \(\mathcal{H}_{\pi_2}^2\) is a generalized graded algebra homomorphism (isomorphism) \(\{f_p\colon H_p^1 \to H_{\varPhi(p)}^2\}_{p\in\pi_1}\) from \(H_{\pi_1}^1\) to \(H_{\pi_2}^2\), satisfying that for any \(p\in\pi_1\), \(f_p\) is a coalgebra homomorphism from \(\mathcal{H}_p^1\) to \(\mathcal{H}_{\varPhi(p)}^2\), and for any \(p\in\pi_1\), \(f_{p^{-1}} \circ S_p^1 = S_{\varPhi(p)}^2 \circ f_p.\)
    
\item When \(\pi_1 = \pi_2 = \pi\), the above cases reduce to the classical \(\pi\)-algebra homomorphism (isomorphism) and Hopf \(\pi\)-algebra homomorphism (isomorphism), respectively.
\end{enumerate}

\item A Hopf \(\pi\)-algebra \(\mathcal{H}_\pi\)  is called commutative if \(\pi\) is an abelian group and for all \(p, q \in \pi\), \(m_{q,p}\circ \sigma_{p,q} = m_{p,q}\); it is called cocommutative if for every \(p \in \pi\), \(\Delta_p = \Delta_p^{\mathrm{cop}}\).

\item A \(\pi\)-subalgebra \(\bigl(\{G_p\}_{p \in \pi}, \{m_{p,q}\}_{p,q \in \pi}, \mu\bigr)\) of the \(\pi\)-algebra \(H_\pi\) satisfies the following conditions:
\begin{enumerate}[label=(\roman*)]
\item For all \(p\in\pi\), \(G_p\) is a subspace of \(H_p\);
\item For all \(p, q \in \pi\), \(m_{p,q}(G_p \otimes G_q) \subseteq G_{pq}\);
\item \(\mu(\mathbf{1}_{K}) \in G_{e}\).
\end{enumerate}

\item A Hopf \(\pi\)-subalgebra 
\(\bigl(\{G_p, \Delta_p, \varepsilon_p\}_{p\in\pi}, \{m_{p,q}\}_{p,q\in\pi}, \mu, \{S_p\}_{p\in\pi}\bigr)\)
of the Hopf \(\pi\)-algebra \(\mathcal{H}_\pi\) satisfies the following conditions:
\begin{enumerate}[label=(\roman*)]
\item For any \(p\in\pi\), \((G_p, \Delta_p, \varepsilon_p)\) is a subcoalgebra of \(\mathcal{H}_p\);
\item \(\bigl(\{G_p\}_{p\in\pi}, \{m_{p,q}\}_{p,q\in\pi}, \mu\bigr)\) is a \(\pi\)-subalgebra of \(H_\pi\);
\item For any \(p\in\pi\), \(S_p(G_p) \subseteq G_{p^{-1}}\).
\end{enumerate}
\end{enumerate}
\end{Definition}

Based on the above definitions, we first present the conclusion on the opposite \(\pi\)-algebra, a fundamental construction in the graded setting.

\begin{Lemma}
If \(H_\pi\) is a \(\pi\)-algebra, then $$H_\pi^{\mathrm{op}}:=\bigl(\{H_p^{\mathrm{op}}= H_{p^{-1}}\}_{p \in \pi}, \{m_{p,q}^{\mathrm{op}} = m_{q^{-1},p^{-1}}\circ\sigma_{p^{-1},q^{-1}}\}_{p,q \in \pi}, \mu\bigr)$$ is also a \(\pi\)-algebra.
\end{Lemma}

\begin{proof}
\begin{align*}
m^{\mathrm{op}}_{pq,r} \circ (m^{\mathrm{op}}_{p,q} \otimes \mathrm{id}_{H_{r^{-1}}}) & = m_{r^{-1},q^{-1}p^{-1}} \circ \sigma_{q^{-1}p^{-1},r^{-1}} \circ (m_{q^{-1},p^{-1}} \circ \sigma_{p^{-1},q^{-1}} \otimes \mathrm{id}_{H_{r^{-1}}}) \\
& = m_{r^{-1}q^{-1},p^{-1}} \circ (m_{r^{-1},q^{-1}} \circ \sigma_{q^{-1},r^{-1}} \otimes \mathrm{id}_{H_{p^{-1}}}) \\
&= m_{r^{-1},\, q^{-1}p^{-1}} \circ \bigl( \mathrm{id}_{H_{r^{-1}}} \otimes m_{q^{-1},\,p^{-1}} \bigr) \circ \sigma_{13} \\
&= m_{r^{-1}q^{-1},\, p^{-1}} \circ \bigl( m_{r^{-1},\,q^{-1}} \otimes \mathrm{id}_{H_{p^{-1}}} \bigr) \circ \sigma_{13}\\
& = m_{r^{-1}q^{-1},p^{-1}} \circ \sigma_{p^{-1},r^{-1}q^{-1}} \circ (\mathrm{id}_{H_{p^{-1}}} \otimes m_{r^{-1},q^{-1}} \circ \sigma_{q^{-1},r^{-1}}) \\
& = m^{\mathrm{op}}_{p,qr} (\mathrm{id}_{H_{p^{-1}}} \otimes m^{\mathrm{op}}_{q,r})
\end{align*}
\begin{align*}
m^{\mathrm{op}}_{p,e} \circ (\mathrm{id}_{H_{p^{-1}}} \otimes \mu) \circ \phi
&= m_{e,p^{-1}} \circ \sigma_{p^{-1},e} \circ (\mathrm{id}_{H_{p^{-1}}} \otimes \mu) \circ \phi \\
&= m_{e,p^{-1}} \circ (\mu \otimes \mathrm{id}_{H_{p^{-1}}}) \circ \phi \\
&= \mathrm{id}_{H_{p^{-1}}} \\
&= m_{p^{-1},e} \circ (\mathrm{id}_{H_{p^{-1}}} \otimes \mu) \circ \phi \\
&= m_{p^{-1},e} \circ \sigma_{e,p^{-1}} \circ (\mu \otimes \mathrm{id}_{H_{p^{-1}}}) \circ \phi
= m^{\mathrm{op}}_{e,p} \circ (\mu \otimes \mathrm{id}_{H_{p^{-1}}}) \circ \phi
\end{align*}
\end{proof}

Extending the conclusion on opposite structures to the case of Hopf \(\pi\)-algebras, we obtain the following result.

\begin{Lemma}
If \(\mathcal{H}_\pi\) is a Hopf \(\pi\)-algebra, then
$
\mathcal{H}_\pi^{\mathrm{op},\mathrm{cop}} := \Bigl(
\bigl\{ \mathcal{H}_p^{\mathrm{op},\mathrm{cop}} = H_{p^{-1}},\ 
\Delta_p^{\mathrm{op},\mathrm{cop}} = \Delta_{p^{-1}}^{\mathrm{cop}},\ 
\varepsilon_p^{\mathrm{op},\mathrm{cop}} =\allowbreak \varepsilon_{p^{-1}} \bigr\}_{p\in\pi},\ \allowbreak 
\bigl\{ m_{p,q}^{\mathrm{op},\mathrm{cop}} = m_{q^{-1},p^{-1}}\allowbreak \circ \sigma_{p^{-1},q^{-1}} \bigr\}_{p,q\in\pi},\ \mu,
\bigl\{ S_p^{\mathrm{op},\mathrm{cop}} = S_{p^{-1}} \bigr\}_{p\in\pi}
\Bigr)
$
is also a Hopf \(\pi\)-algebra.
\end{Lemma}

\begin{proof}
By Lemma 1, \(\Bigl( \bigl\{ H_p^{\mathrm{op},\mathrm{cop}} \bigr\}_{p\in\pi},\ \bigl\{ m_{p,q}^{\mathrm{op},\mathrm{cop}} \bigr\}_{p,q\in\pi},\ \mu \Bigr)\) is a \(\pi\)-algebra, and
\(\Bigl\{ \bigl( H_p^{\mathrm{op},\mathrm{cop}},\ \allowbreak\Delta_p^{\mathrm{op},\mathrm{cop}},\ \allowbreak \varepsilon_p^{\mathrm{op},\mathrm{cop}} \bigr) \Bigr\}_{p\in\pi}\)
is a family of coalgebras. For any \(p,q\in\pi\), we have:
\begin{align*}
\Delta_{pq}^{\mathrm{op},\mathrm{cop}} \circ m_{p,q}^{\mathrm{op},\mathrm{cop}}
&= \Delta^{\mathrm{cop}}_{q^{-1}p^{-1}} \circ m_{q^{-1},p^{-1}} \circ \sigma_{p^{-1},q^{-1}} \\
&= \sigma\circ\Delta_{q^{-1}p^{-1}} \circ m_{q^{-1},p^{-1}} \circ \sigma_{p^{-1},q^{-1}} \\
&= \sigma\circ\bigl( m_{q^{-1},p^{-1}} \otimes m_{q^{-1},p^{-1}} \bigr) \circ \sigma_{23} \circ \bigl( \Delta_{q^{-1}} \otimes \Delta_{p^{-1}} \bigr) \circ \sigma_{p^{-1},q^{-1}} \\
&= \sigma\circ\bigl( m_{q^{-1},p^{-1}} \otimes m_{q^{-1},p^{-1}} \bigr) \circ \bigl( \sigma_{p^{-1},q^{-1}} \otimes \sigma_{p^{-1},q^{-1}} \bigr) \circ \sigma_{23} \circ \bigl( \Delta_{p^{-1}} \otimes \Delta_{q^{-1}} \bigr) \\
&= \bigl( m_{q^{-1},p^{-1}} \circ \sigma_{p^{-1},q^{-1}} \,\otimes\, m_{q^{-1},p^{-1}} \circ \sigma_{p^{-1},q^{-1}} \bigr) \circ \sigma_{23} \circ \bigl( \Delta^{\mathrm{cop}}_{p^{-1}} \otimes \Delta^{\mathrm{cop}}_{q^{-1}} \bigr) \\
&= \bigl( m_{p,q}^{\mathrm{op},\mathrm{cop}} \otimes m_{p,q}^{\mathrm{op},\mathrm{cop}} \bigr) \circ \sigma_{23} \circ \bigl( \Delta_p^{\mathrm{op},\mathrm{cop}} \otimes \Delta_q^{\mathrm{op},\mathrm{cop}} \bigr),
\end{align*}
\begin{align*}
\varepsilon_{pq}^{\mathrm{op},\mathrm{cop}} \circ m_{p,q}^{\mathrm{op},\mathrm{cop}}&= \varepsilon_{q^{-1}p^{-1}} \circ m_{q^{-1},p^{-1}} \circ \sigma_{p^{-1},q^{-1}} = m_\mathbb{K} \circ \bigl( \varepsilon_{q^{-1}} \otimes \varepsilon_{p^{-1}} \bigr) \circ \sigma_{p^{-1},q^{-1}} \\
&= m_\mathbb{K} \circ \bigl( \varepsilon_{p^{-1}} \otimes \varepsilon_{q^{-1}} \bigr) = m_\mathbb{K} \circ \bigl( \varepsilon_p^{\mathrm{op},\mathrm{cop}} \otimes \varepsilon_q^{\mathrm{op},\mathrm{cop}} \bigr),
\end{align*}
meanwhile,
$
\Delta_e^{\mathrm{op},\mathrm{cop}} \circ \mu = (\mu \otimes \mu) \circ \Delta_{\mathbb{K}}$
and
$\varepsilon_e^{\mathrm{op},\mathrm{cop}} \circ \mu = \mathrm{id}_{\mathbb{K}},
$
i.e. both \(\bigl\{ m_{p,q}^{\mathrm{op},\mathrm{cop}} \bigr\}_{p,q\in\pi}\) and \(\mu\) are coalgebra homomorphisms. 

For any \(p\in\pi\), we have:
\begin{align*}
m_{p^{-1},p}^{\mathrm{op},\mathrm{cop}} \circ \bigl( S_p^{\mathrm{op},\mathrm{cop}} \otimes  \mathrm{id}_{H_p} \bigr) \circ \Delta_p^{\mathrm{op},\mathrm{cop}}
&= m_{p^{-1},p} \circ \sigma_{p,p^{-1}} \circ \bigl( S_{p^{-1}} \otimes  \mathrm{id}_{H_p} \bigr) \circ \Delta_{p^{-1}} \\
&= m_{p^{-1},p} \circ \sigma_{p,p^{-1}} \circ \bigl( S_{p^{-1}} \otimes  \mathrm{id}_{H_p} \bigr) \circ \sigma \circ \Delta_{p^{-1}} \\
&= m_{p^{-1},p} \circ \bigl( \mathrm{id}_{H_p} \otimes S_{p^{-1}}  \bigr) \circ \Delta_{p^{-1}} \\
&= u \circ \varepsilon_{p^{-1}} = u \circ \varepsilon_p^{\mathrm{op},\mathrm{cop}},
\end{align*}
similarly,
$
m_{p,p^{-1}}^{\mathrm{op},\mathrm{cop}} \circ \bigl( \mathrm{id}_p \otimes S_p^{\mathrm{op},\mathrm{cop}} \bigr) \circ \Delta_p^{\mathrm{op},\mathrm{cop}}= u \circ \varepsilon_p^{\mathrm{op},\mathrm{cop}}.
$
\end{proof}

The following lemma will be used in the construction of concrete examples in subsequent sections.

\begin{Lemma}
(\cite{radford2012hopf})Let \(U, V\) be two vector spaces. Define a linear operator
 $$\pi: \operatorname{Hom}_{\mathbb{K}}(U, K) \otimes V \to \operatorname{Hom}_{\mathbb{K}}(U, V),\pi(u^* \otimes v)(u) = u^*(u)v$$
 for any \(u^* \otimes v \in \operatorname{Hom}_{\mathbb{K}}(U, K) \otimes V\) and \(u \in U\).  
 
 Then \(\pi\) is a linear isomorphism from \(\operatorname{Hom}_{\mathbb{K}}(U, K) \otimes V\) onto \(\operatorname{Hom}_f(U, V)\),  where \(\operatorname{Hom}_f(U, V)\) denotes the space of finite-rank linear operators from \(U\) to \(V\).
\end{Lemma}

\section{Convolution \(\pi\)-Algebras}

This section establishes the complete theoretical framework of convolution \(\pi\)-algebras, which serves as the core technical tool for the analysis of Hopf \(\pi\)-algebra structures and Rota--Baxter operators throughout the paper.

We present the definition of convolution \(\pi\)-algebra and take it as the main tool to study Hopf \(\pi\)-algebra. 

\begin{Theorem}
Let \(H_\pi\) be a \(\pi\)-algebra and let \(\mathcal{C}\) be a coassociative counital coalgebra. 
Define a family of vector spaces \(\{\operatorname{Hom}_{\mathbb{K}}(C, H_p)\}_{p \in \pi}\) and a family of bilinear operators
\begin{equation*}
\begin{aligned}
\bigl\{\prescript{\mathcal{C}}{H_\pi}{m}_{p,q}^{\#}: {}&\operatorname{Hom}_{\mathbb{K}}(C, H_p) \otimes \operatorname{Hom}_{\mathbb{K}}(C, H_q) \to \operatorname{Hom}_{\mathbb{K}}(C, H_{pq}), \\
&\prescript{\mathcal{C}}{H_\pi}{m}_{p,q}^{\#}(f \otimes g) = f \prescript{\mathcal{C}}{H_{\pi}}\#_{p,q} g = m_{p,q}\circ(f \otimes g)\circ\Delta\bigr\}_{p,q \in \pi}
\end{aligned}
\end{equation*}
and a linear operator $$\prescript{\mathcal{C}}{H_\pi}{\mu}^{\#}: \mathbb{K} \to \operatorname{Hom}_{\mathbb{K}}(C, H_e), \prescript{\mathcal{C}}{H_\pi}{\mu}^{\#}(1_K) = \mu \circ\varepsilon$$
,then \(\bigl(\{\operatorname{Hom}_{\mathbb{K}}(C, H_p)\}_{p \in \pi},\allowbreak \{\prescript{\mathcal{C}}{H_\pi}{m}_{p,q}^{\#}\}_{p,q \in \pi},\allowbreak \prescript{\mathcal{C}}{H_\pi}{\mu}^{\#}\bigr)\) forms a \(\pi\)-algebra. It is called the convolution \(\pi\)-algebra induced by the coalgebra \(\mathcal{C}\) and the \(\pi\)-algebra \(H_\pi\), denoted by \(\operatorname{Conv}(\mathcal{C}, H_\pi)\). This notation is used throughout the rest of the paper.
\end{Theorem}

\begin{proof}
Take any \(p, q, r \in \pi\), \(f \in \operatorname{Hom}_{\mathbb{K}}(C, H_p)\), \(g \in \operatorname{Hom}_{\mathbb{K}}(C, H_q)\), \(j \in \operatorname{Hom}_{\mathbb{K}}(C, H_r)\). Then
\begin{align*}
\bigr(m_{pq,r}^\# \circ (\prescript{\mathcal{C}}{H_\pi}{m}_{p,q}^{\#} \otimes \mathrm{id}_{H_r}) \bigr) (f \otimes g \otimes j)
&= m_{pq,r} \circ \bigl( m_{p,q} \circ (f \otimes g) \circ \Delta \otimes j \bigr) \circ \Delta \\
&= m_{pq,r} \circ (m_{p,q} \otimes \mathrm{id}_{H_r}) \circ (f \otimes g \otimes j) \circ (\Delta \otimes \mathrm{id}_C) \circ \Delta \\
&= m_{p,qr} \circ (\mathrm{id}_{H_p} \otimes m_{q,r}) \circ (f \otimes g \otimes j) \circ (\mathrm{id}_C \otimes \Delta) \circ \Delta \\
&= m_{p,qr} \circ \bigl( f \otimes m_{q,r} \circ (g \otimes j) \circ \Delta \bigr) \circ \Delta\\
&= \bigr(\prescript{\mathcal{C}}{H_\pi}{m}_{p,qr}^{\#} \circ (\mathrm{id}_{H_p} \otimes \prescript{\mathcal{C}}{H_\pi}{m}_{q,r}^{\#}))  (f \otimes g \otimes j),
\end{align*}
\begin{equation*}
\begin{aligned}
\bigl(m_{p,e}^{\#} \circ (\mathrm{id}_{\operatorname{Hom}_{\mathbb{K}}(C,H_p)} \otimes \prescript{\mathcal{C}}{H_\pi}{\mu}^{\#}) \circ \phi\bigr)(f)&= m_{p,e}^{\#}\bigl(f \otimes \mu\circ\varepsilon\bigr) = m_{p,e} \circ (f \otimes \mu\circ\varepsilon) \circ \Delta\\
&=m_{p,e} \circ \bigl(\mu\otimes \operatorname{id}_{\operatorname{Hom}_{\mathbb{K}}(C, H_p)} \bigr) \circ \phi \circ f \circ \phi^{-1} \circ (\varepsilon\otimes \operatorname{id}_C) \circ \Delta\\
&= \operatorname{id}_{\operatorname{Hom}_{\mathbb{K}}(C, H_p)} (f)\\
&=m_{p,e} \circ \bigl( \operatorname{id}_{\operatorname{Hom}_{\mathbb{K}}(C, H_p)} \otimes \mu \bigr) \circ \phi \circ f \circ \phi^{-1} \circ (\operatorname{id}_C \otimes \varepsilon) \circ \Delta\\
&= m_{p,e} \circ (\mu\circ\varepsilon \otimes f) \circ \Delta = m_{p,e}^{\#}\bigl(\mu\circ\varepsilon \otimes f\bigr) \\
&= \bigl(m_{p,e}^{\#} \circ (\prescript{\mathcal{C}}{H_\pi}{\mu}^{\#} \otimes \mathrm{id}_{\operatorname{Hom}_{\mathbb{K}}(C,H_p)}) \circ \phi\bigr)(f).
\end{aligned}
\end{equation*}		
\end{proof}
	
\begin{Example}
Let \(H_\pi\) be a \(\pi\)-algebra.
\begin{enumerate}[label=(\alph*)]	
\item Given any coassociative  counital coalgebra \(\mathcal{C}\), define a family of linear operators
	$$\{\beta_p : \operatorname{Hom}_{\mathbb{K}}(C, K) \otimes H_p \to \operatorname{Hom}_{\mathbb{K}}(C, H_p),\allowbreak  \beta_p(c^* \otimes h)(c) = c^*(c)h \}_{p \in \pi},$$
	then \((\{\operatorname{Im}(\beta_p)\}_{p \in \pi},\allowbreak \prescript{\mathcal{C}}{H_\pi}{m}_{p,q}^{\#},\allowbreak \prescript{\mathcal{C}}{H_\pi}{\mu}^{\#})\) is a \(\pi\)-subalgebra of \(\operatorname{Conv}(\mathcal{C}, H_\pi)\).
	
\item In particular, take \(\mathcal{C} = \mathbb{K}\). Then there is an isomorphism of \(\pi\)-algebras between \(\operatorname{Conv}(\mathbb{K}, H_\pi)\) and \(H_\pi\).
\end{enumerate}
\end{Example}

\begin{proof}
\begin{enumerate}[label=(\alph*)]	
\item By Lemma 3, \(\operatorname{Im}(\beta_p)\) is a linear subspace of \(\operatorname{Hom}_{\mathbb{K}}(C, H_p)\) for each \(p\in\pi\). 
Take arbitrary \(c^*, d^*\in C^*\), \(h\in H_p\), \(g\in H_q\) and \(n\in C\). We have:
\[
\begin{aligned}
{}_{H_\pi}^C m_{p,q}^\# \big( \beta_p(c^* \otimes h) \otimes \beta_q(d^* \otimes g) \big) (n)
&= \bigr(m_{p,q} \circ \big( \beta_p(c^* \otimes h) \otimes \beta_q(d^* \otimes g) \big) \circ \Delta\bigr) (n) \\
&= c^*(n_{(1)})\, d^*(n_{(2)})\, h g = (c^* d^*)(n)\, h g \\
&= \beta_{pq}\big( c^* d^* \otimes h g \big) (n) \in \operatorname{Im}(\beta_{pq}),
\end{aligned}
\]
and $
{}_{H_\pi}^C \mu^\#(1_\mathbb{K}) = \mu \circ \varepsilon = \beta_e\big(\varepsilon \otimes \mu(1_\mathbb{K})\big) \in \operatorname{Im}(\beta_e).
$

\item Take a family of linear space isomorphisms $$\{\psi_p \colon \operatorname{Hom}_{\mathbb{K}}(\mathbb{K}, H_p) \to H_p\}_{p\in\pi},\psi_p(f) = f(1_{\mathbb{K}}).$$
For any \(f \in \operatorname{Hom}_{\mathbb{K}}(\mathbb{K}, H_p)\) and \(g \in \operatorname{Hom}_{\mathbb{K}}(\mathbb{K}, H_q)\), it holds that
\[
\begin{aligned}
\big( m_{p,q} \circ (\psi_p \otimes \psi_q) \big) (f \otimes g)
&= m_{p,q}\big( f(1_{\mathbb{K}}) \otimes g(1_{\mathbb{K}}) \big) \\
&= \psi_{pq}\big( m_{p,q} \circ (f \otimes g) \circ \Delta_{\mathbb{K}} \big) = \big( \psi_{pq} \circ {}_{H_{\pi}}^{\mathbb{K}} m_{p,q}^{\#} \big) (f \otimes g)
\end{aligned}
\]
and \(\psi_e \circ {}_{H_\pi}^{\mathbb{K}} \mu^\# = (\mu \circ \varepsilon_{\mathbb{K}})(1_{\mathbb{K}}) = \mu.\)
\end{enumerate}
\end{proof}

In particular, the antipode axiom of Hopf \(\pi\)-algebras can be equivalently characterized by means of convolution inverses.

\begin{Remark}
Let \(\mathcal{H}_\pi\) be a Hopf \(\pi\)-algebra.
\begin{enumerate}[label=(\alph*)]	
\item Then for each \(p\in\pi\), we have a coassociative counital coalgebra \(\mathcal{H}_p\) and a \(\pi\)-algebra \(H_\pi\). Consequently, we obtain the convolution \(\pi\)-algebra \(\operatorname{Conv}(\mathcal{H}_p,H_\pi )\). 
\item It is easy to see that the antipode compatibility condition of the Hopf \(\pi\)-algebra is equivalent to the following:
 For each \(p \in \pi\), \(S_p\) and \(\mathrm{id}_p\) are convolution inverses of each other in \(\operatorname{Conv}(\mathcal{H}_p,H_\pi )\).
\end{enumerate}
\end{Remark}

Homomorphisms between convolution \(\pi\)-algebras are discussed next, which can be naturally induced by homomorphisms of the underlying structures.

\begin{Theorem}
Let \(\alpha:\mathcal{C} \to \mathcal{D}\) be a coalgebra homomorphism, 
and let \(\{\beta_p \colon H_p^1 \to H_{\varPhi(p)}^2\}_{p\in \pi_1}\) be a generalized graded algebra homomorphism from \(H_{\pi_1}^1\) to \(H_{\pi_2}^2\). Then
\begin{enumerate}[label=(\alph*)]	
\item $$\{\alpha_p^\# : \operatorname{Hom}_{\mathbb{K}}(D, H_p^1) \to \operatorname{Hom}_{\mathbb{K}}(C, H_p^1),\ \alpha_p^\#(f) = f \circ \alpha\}_{p \in \pi_1}$$ 
	is a  \(\pi_1\)-algebra homomorphism from 
	\(\operatorname{Conv}(\mathcal{D},H_{\pi_1}^1 )\)
	to 
\(\operatorname{Conv}(\mathcal{C},H_{\pi_1}^1 )\).
	
\item 
	$$\{\beta_p^\# : \operatorname{Hom}_{\mathbb{K}}(C, H_p^1) \to \operatorname{Hom}_{\mathbb{K}}(C, H_{\varPhi (p)}^2),\ \beta_p^\#(f) = \beta_p \circ f \}_{p\in\pi_1}$$ 
	is a generalized graded algebra homomorphism from 
	\(\operatorname{Conv}(\mathcal{C},H_{\pi_1}^1) \) 
	to 
	\(\operatorname{Conv}(\mathcal{C},H_{\pi_2}^2) \).
\end{enumerate}
\end{Theorem}

\begin{proof}
(1) For any \(p, q \in \pi_1\), \(f \in \operatorname{Hom}_{\mathbb{K}}(D,H_p^1)\), \(g \in \operatorname{Hom}_{\mathbb{K}}(D, H_q^1)\), we have
\begin{align*}
\left( \prescript{\mathcal{C}}{H_{\pi_1}^1}{m}_{p,q}^{\#} \circ (\alpha_p^{\#} \otimes \alpha_q^{\#}) \right)(f \otimes g) 
&= \prescript{\mathcal{C}}{H_{\pi_1}^1}{m}_{p,q}^{\#}\bigr((f \otimes g)\circ(\alpha \otimes \alpha)\bigr)\\
&= m_{p,q}^{1}\circ(f \otimes g)\circ(\alpha \otimes \alpha)\circ\triangle_C \\
&= m_{p,q}^{1}\circ(f \otimes g)\circ \triangle_C\circ\alpha = \bigr(\prescript{\mathcal{D}}{H_{\pi_1}^1}{m}_{p,q}^{\#}(f \otimes g)\bigr)\circ\alpha \\
&= \alpha_{p,q}^{\#} \bigr(\prescript{\mathcal{D}}{H_{\pi_1}^1}{m}_{p,q}^{\#}(f \otimes g)\bigr) 
= \left( \alpha_{p,q}^{\#} \circ \prescript{\mathcal{C}}{H_{\pi_1}^1}{m}_{p,q}^{\#} \right)(f \otimes g)
\end{align*}
and \(\alpha_e^{\#}(\mu^1 \circ\varepsilon_D) = \mu^1 \circ\varepsilon_D\circ \alpha = \mu^A\circ \varepsilon\).

(2) For any \(p, q \in \pi_1\), \(f \in \operatorname{Hom}_{\mathbb{K}}(C, H_p^1)\), \(g \in \operatorname{Hom}_{\mathbb{K}}(C, H_q^1)\), we have
\begin{equation*}
\begin{aligned}
\phantom{=} \bigl( \beta_{pq}^\# \circ \prescript{\mathcal{C}}{H_{\pi_1}^1}{m}_{p,q}^{\#} \bigr) (f \otimes g)
&= \beta_{pq}^\# \bigl( \prescript{\mathcal{C}}{H_{\pi_1}^1}{m}_{p,q}^{\#} (f \otimes g) \bigr) \\
&= \beta_{pq} \circ m_{p,q}^1 \circ (f \otimes g) \circ \Delta\\
&= m_{\varPhi(p),\varPhi(q)}^2 \circ (\beta_p \otimes \beta_q) \circ (f \otimes g) \circ \Delta \\
&= \prescript{\mathcal{C}}{H_{\pi_2}^2}{m}_{\varPhi(p),\varPhi(q)}^{\#} \bigl( (\beta_p \otimes \beta_q) \circ (f \otimes g) \bigr) \\
&=\bigl( \prescript{\mathcal{C}}{H_{\pi_2}^2}{m}_{\varPhi(p),\varPhi(q)}^{\#} \circ (\beta_p^\# \otimes \beta_q^\#) \bigr) (f \otimes g) 
\end{aligned}
\end{equation*}
and \(\beta_e^{\#}(\mu^A \circ\varepsilon) = \beta_e \circ\mu^A \circ\varepsilon = \mu^B \circ\varepsilon\).
\end{proof}

Based on the above results on homomorphism induction, the explicit expressions of convolution inverses for two typical classes of maps are presented.

\begin{Theorem}
Suppose \(\mathcal{H}_{\pi_1}\) is a Hopf \(\pi_1\)-algebra.
\begin{itemize}
\item[(1)]  If \(A_{\pi_2}\) is a \(\pi_2\)- algebra , and \(\{f_p : H_p \to A_{\varPhi(p)}\}_{p \in {\pi_1}}\) is a generalized graded algebra homomorphism from \(H_{\pi_1}\) to \(A_{\pi_2}\),then for any \(p\in\pi_1\), \(f_p\) has a convolution inverse \(f_{p^{-1}} \circ S_p\) in \(\operatorname{Conv}(\mathcal{H}_p,A_{\pi_2} )\).

\item[(2)]  If \(\mathcal{C}\) is a coalgebra and \(f_p : C \to H_p\) is a coalgebra homomorphism from \(\mathcal{C}\) to \(\mathcal{H}_p\), then \(f_p \in \operatorname{Hom}_{\mathbb{K}}(C, H_p)\) has a convolution inverse \(S_p \circ f_p\) in \(\operatorname{Conv}(\mathcal{C},H_{\pi_1} )\).
\end{itemize}
\end{Theorem}
\begin{proof}
(1) For any \(p \in \pi_1\), we have the coassociative counital coalgebra \(\mathcal{H}_p\). Since \(\{f_p: H_p \to A_{\varPhi(p)}\}_{p \in \pi_1}\) is a generalized algebra homomorphism from \(H_{\pi_1}\) to \(A_{\pi_2}\), by Theorem 2, there is a generalized graded algebra homomorphism \(\{f_p^{\#}\}_{p \in \pi_1}\) from \(\operatorname{Conv}(\mathcal{H}_p,H_{\pi_1} )\) to \(\operatorname{Conv}(\mathcal{H}_p,A_{\pi_2} )\).

Since \(S_p\) and \(\mathrm{id}_p\) are convolution inverses of each other in  \(\operatorname{Conv}(\mathcal{H}_p,H_{\pi_1} )\), it follows that \(f_{p^{-1}}^{\#}(S_p) = f_{p^{-1}}\circ S_p\) and \(f_p^{\#} \mathrm({id}_p) = f_p\) are convolution inverses of each other in \(\operatorname{Conv}(\mathcal{H}_p,A_{\pi_2} )\).

(2) It is known that \(f_p\) is a coalgebra homomorphism from \(\mathcal{C}\) to \(\mathcal{H}_p\), and \(H_{\pi_1}\) is a \(\pi\)-algebra, by Theorem 2, there is a convolution \(\pi_1\)-algebra homomorphism \(\{f_{r}^{p,\#}\}_{r \in \pi}\) from \(\operatorname{Conv}(\mathcal{H}_p,H_{\pi_1} )\) to \(\operatorname{Conv}(\mathcal{C},H_{\pi_1} )\).

Since \(S_p\) and \(id_p\) are convolution inverses of each other in \(\operatorname{Conv}(\mathcal{H}_p,H_{\pi} )\), it follows that \(f_{p^{-1}}^{\,p,\mathbin{\#}}(S_p)= S_p \circ f_p\) and \(f_{p^{-1}}^{\,p,\mathbin{\#}}(id_p) = f_p\) are convolution inverses of each other in \(\operatorname{Conv}(\mathcal{C},H_{\pi} )\).
\end{proof}

As a direct corollary, the criterion for the antipode compatibility of Hopf \(\pi\)-algebra homomorphisms can be derived.

\begin{Corollary}
Let \(\mathcal{H}^1_{\pi_1}\) be a Hopf \(\pi_1\)-algebra and \(\mathcal{H}^2_{\pi_2}\) be a Hopf \(\pi_2\)-algebra, and let \(\{f_p: H_p^1 \to H^2_{\Phi(p)}\}_{p\in\pi_1}\) be a generalized graded algebra homomorphism from \(\mathcal{H}^1_{\pi_1}\) to \(\mathcal{H}^2_{\pi_2}\). If for each \(p\in\pi_1\), \(f_p\) is a coalgebra homomorphism from \(\mathcal{H}_p^1\) to \(\mathcal{H}^2_{\Phi(p)}\), then for any \(p\in\pi_1\), \(f_{p^{-1}} \circ S_p^1 = S^2_{\Phi(p)} \circ f_p.\)
\end{Corollary}

\begin{proof}
For any \(p \in \pi_1\), \(\{f_p\}_{p \in \pi_1}\) is a generalized graded homomorphism from \(H^1_{\pi_1}\) to \(H_{\pi_2}^2\). By Theorem 3, \(f_p : H^1_p \to H^2_{\varPhi(p)}\) is convolution inverse to \(f_{p^{-1}}\circ S^1_p\) in  \(\operatorname{Conv}(\mathcal{H}^1_p,H_{\pi}^2 )\).

For any \(p \in \pi\), \(f_p\) is a coalgebra homomorphism from \(\mathcal{H}^1_p\) to \(\mathcal{H}_{\varPhi(p)}^2\). By Theorem 3, \(f_p : H^1_p \to H_{\varPhi(p)}^2\) is convolution inverse to \(S_{\varPhi(p)}^2 \circ f_p\) in \(\operatorname{Conv}(\mathcal{H}^1_p,H_{\pi_2}^2 )\).

In summary, for any \(p \in \pi_1\), \(f_{p^{-1}} \circ S_p^1 = S_{\varPhi(p)}^2\circ f_p\).
\end{proof}

The anti-algebra homomorphism property of convolution inverses is stated in the following lemma.

\begin{Lemma}
Let \(\mathcal{H}_{\pi_1}\) be a Hopf \(\pi_1\)-algebra and  \(A_{\pi_2}\) be a \(\pi_2\)-algebra. If \(\bigl\{f_p: H_p \to A_{\varPhi(p)}\bigr\}_{p\in\pi_1}\) is a generalized graded algebra homomorphism
from \(H_{\pi_1}\) to \(A_{\pi_2}\) and for each \(p \in \pi_1\), \(f_p: H_p \to A_{\varPhi(p)}\) has a convolution inverse \(f_p^{-1}: H_p \to A_{\varPhi(p)}^{op}\) in \(\operatorname{Conv}(\mathcal{H}_p,A_{\pi_2} )\); then \(\{f_p^{-1}\}_{p\in\pi_1}\) is a generalized graded algebra homomorphism from \(H_{\pi_1}\) to \(A_{\pi_2}^{op}\).
\end{Lemma}
	
\begin{proof}
For any \(p,q \in \pi\), the tensor product coalgebra of \(\mathcal{H}_p\) and \(\mathcal{H}_q\) is 
$$\mathcal{H}_p\otimes \mathcal{H}_q:=(H_p \otimes H_q, \Delta_{p,q}^\otimes=\sigma_{23}\circ (\Delta_p\otimes\Delta_q), \varepsilon_{p,q}^\otimes=m_{\mathbb{K}}\circ (\varepsilon_p\otimes\varepsilon_q)).$$ 
Consequently, we have \(\operatorname{Conv}(\mathcal{H}_p\otimes \mathcal{H}_q,A_{\pi_2} )\). For any \(p,q\in \pi\), 
\begin{align*}
&\prescript{\mathcal{H}_p\otimes\mathcal{H}_q}{A_\pi}m^{\#}_{\varPhi(pq)^{-1},\varPhi(pq)}(f_{pq}^{-1}\circ m_{p,q}\otimes f_{pq}\circ m_{p,q})\\
&=m^{A}_{\varPhi(pq)^{-1},\varPhi(pq)}\circ(f_{pq}^{-1}\circ m_{p,q}\otimes f_{pq}\circ m_{p,q})\circ \Delta_{p,q}^{\otimes}\\
&=m^{A}_{\varPhi(pq)^{-1},\varPhi(pq)}\circ(f_{pq}^{-1}\otimes f_{pq})\circ (m_{p,q}\otimes m_{p,q})\circ \sigma_{23}\circ (\Delta_p\otimes\Delta_q)\\
&=m^{A}_{\varPhi(pq)^{-1},\varPhi(pq)}\circ(f_{pq}^{-1}\otimes f_{pq})\circ \Delta_{p,q}\circ m_{p,q}=\mu^A\circ \varepsilon_{pq}\\
&=\mu\circ m_{\mathbb{K}}\circ (\varepsilon_p\otimes \varepsilon_q)=\mu^A\circ\varepsilon_{p,q}^{\otimes},
\end{align*}
\begin{align*}
&\prescript{\mathcal{H}_p\otimes\mathcal{H}_q}{A_\pi} m_{\varPhi(pq),\varPhi(pq)^{-1}}^{\#}(f_{pq}\circ m_{p,q}\otimes m_{\varPhi(p),\varPhi(q)}^{A,op}\circ (f_p^{-1}\otimes f_q^{-1}))\\
&=m_{\varPhi(pq),\varPhi(pq)^{-1}}^{A}\circ (f_{pq}\circ m_{p,q}\otimes m_{\varPhi(p),\varPhi(q)}^{A,op}\circ (f_p^{-1}\otimes f_q^{-1}))\circ \Delta^{\otimes}_{p,q}\\
&=m_{\varPhi(pq),\varPhi(pq)^{-1}}^{A}\circ (m_{\varPhi(p),\varPhi(q)}^{A}\circ (f_p\otimes f_q)\otimes m_{\varPhi(p),\varPhi(q)}^{A,op}\circ (f_p^{-1}\otimes f_q^{-1}))\circ \Delta^{\otimes}_{p,q}\\
&=m_{\varPhi(pq),\varPhi(pq)^{-1}}^{A}\circ (m_{\varPhi(p),\varPhi(q)}^{A}\otimes m_{\varPhi(p),\varPhi(q)}^{A,op})\circ (f_p\otimes f_q\otimes f_p^{-1}\otimes f_q^{-1})\circ \Delta^{\otimes}_{p,q}\\
&=\mu^A\circ m_{\mathbb{K}}\circ (\varepsilon_p\otimes\varepsilon_q)=\mu^A\circ \varepsilon_{p,q}^{\otimes},\\
\end{align*}
that is, \(f_{pq}^{-1} \circ m_{p,q}\) and \(m_{\varPhi(p)^{-1},\varPhi(q)^{-1}}^{A,op}\circ (f_p^{-1} \otimes f_q^{-1})\) are respectively a convolution left inverse and a convolution right inverse of \(f_{pq}\circ m_{p,q}\) in  \(\operatorname{Conv}(\mathcal{H}_p\otimes \mathcal{H}_q,A_{\pi_2} )\). 

Therefore, for any \(p,q \in \pi\), we have \(m_{\varPhi(p)^{-1},\varPhi(q)^{-1}}^{A,op}\circ(f_p^{-1} \otimes f_q^{-1}) = f_{pq}^{-1}\circ m_{p,q}\).

Since \(f_e\) and \(f_e^{-1}\) are convolution inverses in \(\mathrm{Conv}(H_e, A_{\pi_2})\), it follows that
\begin{equation*}
\begin{aligned}
  \mu^A =& \mu^A \circ \operatorname{id}_{\mathbb{K}} = \mu^A \circ \varepsilon_e \circ \mu = \prescript{H_e}{A_{\pi_2}}m_{e,e}^\# \circ (f_e \otimes f_e^{-1}) \circ \mu =m_{e,e}^A \circ \left( f_e \otimes f_e^{-1}\right) \circ \Delta_e \circ \mu \\
= & m_{e,e}^A \circ \left( f_e \otimes f_e^{-1} \right) \circ \left( \mu \otimes \mu \right) \circ \Delta_{\mathbb{K}} =  m_{e,e}^A \circ \left( f_e \circ \mu \otimes f_e^{-1} \circ \mu \right) \circ \Delta_{\mathbb{K}} \\
= & m_{e,e}^A \circ \left( \mu^A \otimes f_e^{-1} \circ \mu \right) \circ \Delta_{\mathbb{K}} =  m_{e,e}^A \circ \left( \mu^A \circ \operatorname{id}_{\mathbb{K}} \otimes \operatorname{id}_{H_e} \circ f_e^{-1} \circ \mu \right) \circ \Delta_{\mathbb{K}} \\
= & m_{e,e}^A \circ \left( \mu^A \otimes \operatorname{id}_{A_e} \right) \circ \left( \operatorname{id}_{\mathbb{K}} \otimes f_e^{-1} \circ \mu \right) \circ \Delta_{\mathbb{K}} 
=  \operatorname{id}_{A_e} \circ f_e^{-1} \circ \mu =  f_e^{-1} \circ \mu.
\end{aligned}
\end{equation*}
\end{proof}

Correspondingly, the anti-coalgebra homomorphism property of convolution inverses is presented in the following lemma.

\begin{Lemma}
Suppose \(\mathcal{C}\) is a coassociative counital coalgebra, and \(\mathcal{H}_\pi\) is a Hopf \(\pi\)-algebra. If \(f_p \in \operatorname{Hom}_{\mathbb{K}}(C, H_p)\) is a coalgebra homomorphism from \(\mathcal{C}\) to \(\mathcal{H}_p\) and has a convolution inverse \(f_p^{-1} \in \operatorname{Hom}_{\mathbb{K}}(C, H_{p^{-1}})\) in \(\operatorname{Conv}(\mathcal{C},H_{\pi} )\), then \(f_p^{-1}\) is a coalgebra homomorphism from \(\mathcal{C}^{cop}\) to \(\mathcal{H}_{p^{-1}}\).
\end{Lemma}
	
\begin{proof}
It is known that \(H_\pi\) is a \(\pi\)-algebra.
one can easily verify that 
 $$H_\pi\otimes H_\pi:=(\{H_p \otimes H_p\}_{p \in \pi}, \{m_{p,q}^{\otimes}=(m_{p,q}\otimes m_{p,q})\circ \sigma_{23}\}_{p,q \in \pi}, \mu^{\otimes}=(\mu\otimes \mu)\circ\Delta_{\mathbb{K}})$$
  forms a \(\pi\)-algebra, and that \(\{\Delta_p : H_p \to H_p \otimes H_p\}_{p \in \pi}\) is a \(\pi\)-algebra homomorphism from \(H_\pi\) to \(H_{\pi}\otimes H_{\pi}\). 
  
  Consequently, we obtain \(\operatorname{Conv}({\mathcal{C}},H_\pi\otimes H_\pi)\). For any \(p \in \pi\), \(f_p \in \operatorname{Hom}_{\mathbb{K}}(C, H_p)\), we have
\begin{align*}
	\prescript{\mathcal{C}}{H_\pi\otimes H_\pi}m^{\#}_{p,p^{-1}}(\Delta_p\circ f_p\otimes \Delta_{p^{-1}}\circ f_p^{-1})
	&= m_{p,p^{-1}}^{\otimes}\circ (\Delta_p\circ f_p\otimes \Delta_{p^{-1}}\circ f_p^{-1})\circ \Delta \\
	&= (m_{p,p^{-1}}\otimes m_{p,q})\circ \sigma_{23}\circ (\Delta_p\otimes \Delta_{p^{-1}})\circ (f_p\otimes f_p^{-1})\circ \Delta \\
	&= \Delta_e\circ m_{p,q}\circ (f_p\otimes f_p^{-1})\circ \Delta=(\mu\otimes \mu)\Delta_{\mathbb{K}}=\mu^{\otimes}\circ \varepsilon\\
\end{align*}
\begin{equation*}
\begin{aligned}
&\left( {}_{H_\pi}^{c} m_{p^{-1},p}^{\#} \circ \big( (f_p^{-1} \otimes f_p^{-1}) \circ \Delta^{\text{cop}} \otimes \Delta_p \circ f_p \big) \right) (d) \\
&= \left( m_{p^{-1},p}^{\otimes} \circ \big( (f_p^{-1} \otimes f_p^{-1}) \circ \Delta^{\text{cop}} \otimes \Delta_p \circ f_p \big) \circ \Delta \right) (d) \\
&= \left( m_{p^{-1},p}^{\otimes} \circ \big( (f_p^{-1} \otimes f_p^{-1}) \circ \Delta^{\text{cop}} \otimes (f_p \otimes f_p) \circ \Delta \big) \circ \Delta \right) (d) \\
&= \left( (m_{p,p^{-1}} \otimes m_{p,p^{-1}}) \circ \sigma_{23} \circ (f_p^{-1} \otimes f_p^{-1} \otimes f_p \otimes f_p) \circ (\Delta^{\text{cop}} \otimes \Delta) \circ \Delta \right) (d) \\
&= \big( m_{p,p^{-1}} \otimes m_{p,p^{-1}} \big) \Big( f_p^{-1}(d_{(12)}) \otimes f_p(d_{(21)}) \otimes f_p^{-1}(d_{(11)}) \otimes f_p(d_{(22)}) \Big) \\
&= m_{p,p^{-1}}\big( f_p^{-1}(d_{(21)}) \otimes f_p(d_{(22)}) \big) \otimes m_{p,p^{-1}}\big( f_p^{-1}(d_{(1)}) \otimes f_p(d_{(3)}) \big) \\
&= {}_{H_\pi}^{c} m_{p,p^{-1}}^{\#} \big( f_p^{-1} \otimes f_p \big) (d_2) \otimes m_{p,p^{-1}}\big( f_p^{-1}(d_1) \otimes f_p(d_3) \big) \\
&= (\mu \circ \varepsilon)(d_2) \otimes m_{p,p^{-1}}\big( f_p^{-1}(d_1) \otimes f_p(d_3) \big) = 1_H \otimes m_{p,p^{-1}}\big( f_p^{-1}(d_1) \otimes f_p(d_2) \big) \\
&= 1_H \otimes {}_{H_\pi}^{c} m_{p,p^{-1}}^{\#} \big( f_p^{-1} \otimes f_p \big) (d) = (1_H \otimes 1_H) \varepsilon(d) = (\mu^{\otimes} \circ \varepsilon)(d)
\end{aligned}
\end{equation*}
That is, \(\Delta_{p^{-1}}\circ f_p^{-1}\) and \((f_p^{-1} \otimes f_p^{-1})\circ \Delta_{C}^{cop}\) are respectively a convolution right inverse and a convolution left inverse of \(\Delta_p\circ f_p\) in \(\operatorname{Conv}({\mathcal{C}},H_\pi\otimes H_\pi)\). 
 
 Therefore, for any \(p \in \pi, d \in\ C\), we have \(\Delta_{p^{-1}}\circ (f_p^{-1}) = (f_p^{-1} \otimes f_p^{-1})\circ \Delta_{C}^{cop}\) and $$\varepsilon_{p^{-1}} f_p^{-1}(d) = \varepsilon_p f_p(d_1) \varepsilon_{p^{-1}} f_p^{-1}(d_2) = \varepsilon_e (f_p\#_{p,p^{-1}} f_p^{-1})(d) = \varepsilon_e (\mu \varepsilon(d)) = \varepsilon(d).$$
\end{proof}

The following three results constitute the core content of this paper and serve as the cornerstone for constructing examples of Rota-Baxter Hopf \(\pi\)-algebras in the next section.
  
\begin{Theorem}
If \(\mathcal{H}_\pi\) is a Hopf \(\pi\)-algebra, then for any \(p, q \in \pi\), $$S_{pq}\circ m_{pq} = m_{p,q}^{op}\circ (S_p \otimes S_q),\ \ S_e\circ \mu=\mu$$ and for any \(p \in \pi\), $$\Delta_{p^{-1}}\circ S_p = (S_p \otimes S_p)\circ\Delta_{p}^{cop},\ \varepsilon_{p^{-1}}\circ S_p = \varepsilon_p.$$
\end{Theorem}
	
\begin{proof}
It is known that \(\mathcal{H}_\pi\) is a Hopf \(\pi\)-algebra, then for any \(p \in \pi\), \(id_{H_p}\) and \(S_p\) are convolution inverses of each other in \(\operatorname{Conv}(\mathcal{H}_p,H_{\pi} )\); 

Since \(\{id_{H_p} : H_p \to H_p\}_{p \in \pi}\) is an endomorphism of \(H_\pi\) onto itself, by Lemma 4, \(\{S_p\}_{p \in \pi}\) is a \(\pi\)-algebra endomorphism of \(H_\pi\) to \(H_\pi^{op}\).

Since for any \(p \in \pi\), \(id_{H_p}\) is an endomorphism of \(\mathcal{H}_p\) onto itself, by Lemma 5, \(S_p\) is a coalgebra homomorphism from \(\mathcal{H}_p^{cop}\) to \(\mathcal{H}_{p^{-1}}\).
\end{proof}
	
\begin{Theorem}
Suppose \(\mathcal{H}_\pi\) is a Hopf \(\pi\)- algebra.
\begin{itemize}
\item[(1)]  If for any \(p \in \pi\), \(S_p : H_p \to H_{p^{-1}}\) is bijective, then $$(\{H_p, \Delta_p^{cop}, \varepsilon_p\}_{p \in \pi}, \{m_{p,q}\}_{p,q \in \pi}, \mu, \\ \big\{S_p^{cop} = S_{p^{-1}}^{-1}\big\}_{p\in\pi})$$ forms a Hopf \(\pi\)-algebra.
\item[(2)]  If there exists a family of linear operators \(\{S_p^{cop}\}_{p \in \pi}\) such that $$(\{H_p, \Delta_p^{cop}, \varepsilon_p\}_{p \in \pi}, \{m_{p,q}\}_{p,q \in \pi},\\ \mu, \{S_p^{cop}\}_{p \in \pi})$$ forms a Hopf \(\pi\)- algebra, then for any \(p \in \pi\), \(S_p\) is bijective and \(S_p^{cop} = S_{p^{-1}}^{-1}\).
\end{itemize}
\end{Theorem}
	
\begin{proof}
(1) By Theorem 4, \(\{S_p\}_{p \in \pi}\) is a \(\pi\)-algebra isomorphism from \(H_\pi\) to \(H_\pi^{op}\), 
then \(\{F_p = S_p^{-1}\}_{p \in \pi}\) is a \(\pi\)-algebra isomorphism from  \(H_\pi^{op}\) to \(H_\pi\). That is, for any \(p,q \in \pi\), 

\[
m_{p,q}\circ(S_{p^{-1}}\otimes S_{q^{-1}}) =m_{p,q}\circ(F_p\otimes F_q) = F_{qp}\circ m_{p,q}^{\mathrm{op}}= S_{(pq)^{-1}}\circ m_{p,q}^{\mathrm{op}},
\] 
thus for any \(p \in \pi\), we have
\begin{align*}	
m_{p^{-1},p} \circ (S_p^{\mathrm{op}} \otimes \mathrm{id}_{H_p}) \circ \Delta_p^{\mathrm{cop}} 
=&\ m_{p^{-1},p} \circ (S_{p^{-1}}^{-1} \otimes \mathrm{id}_{H_p}) \circ \Delta_p^{\mathrm{cop}} \\
=&\ m_{p^{-1},p} \circ (S_{p^{-1}}^{-1} \otimes S_p^{-1}) \circ (\mathrm{id}_{H_p} \otimes S_p) \circ \Delta_p^{\mathrm{cop}} \\
=&\ S_e^{-1} \circ m_{p^{-1},p}^{\mathrm{op}} \circ (\mathrm{id}_{H_p} \otimes S_p) \circ \Delta_p^{\mathrm{cop}} \\
=&\ S_e \circ m_{p^{-1},p}^{\mathrm{op}} \circ (\mathrm{id}_{H_p} \otimes S_p) \circ \Delta_p^{\mathrm{cop}} \\
=&\ S_e \circ m_{p^{-1},p} \circ (S_p \otimes \mathrm{id}_{H_p}) \circ \Delta_p =\ S_e \circ \mu \circ \varepsilon = \mu \varepsilon,
\end{align*}
similarly, we can prove that
\begin{equation*}
m_{p,p^{-1}} \circ (\mathrm{id}_{H_p} \otimes S_p^{\mathrm{op}}) \circ \Delta_p^{\mathrm{cop}} = \mu \circ \varepsilon.
\end{equation*}
(2) Suppose \((\{H_p, \Delta_p^{cop}, \varepsilon_p\}_{p \in \pi}, \{m_{p,q}\}_{p,q \in \pi}, \mu, \{S_p^{cop}\}_{p \in \pi})\) is a Hopf \(\pi\)- algebra. Then for any \(p \in \pi\), 
\begin{equation*}
m_{p^{-1},p} \circ (S_p^{\mathrm{op}} \otimes \mathrm{id}_{H_p}) \circ \Delta_p^{\mathrm{cop}} = \mu \circ \varepsilon_p = m_{p,p^{-1}} \circ (\mathrm{id}_{H_p} \otimes S_p^{\mathrm{op}}) \circ \Delta_p^{\mathrm{cop}},
\end{equation*}
furthermore, for any \(p \in \pi\), 
\begin{equation*}
m_{p,p^{-1}} \circ (S_{p^{-1}} \circ S_p^{\mathrm{cop}} \otimes S_p) \circ \Delta_p = S_e \circ m_{p,p^{-1}} \circ (\mathrm{id}_{H_p} \otimes S_p^{\mathrm{cop}}) \circ \Delta_p^{\mathrm{cop}} = S_e \circ \mu \circ \varepsilon_p = \mu \circ \varepsilon_p,
\end{equation*}
\(S_{p^{-1}}\circ S_p^{cop}\) is a left inverse of \(S_p\) in \((\{\operatorname{Hom}_{\mathbb{K}}(H_p, H_r)\}_{r \in \pi}, \{m_{r,s}^{p,\#}\}_{r,s \in \pi}, \mu^{p,\#})\), 
therefore \(S_{p^{-1}}\circ S_p^{cop} = \mathrm{id}_{H_p}\).

Similarly, for any \(p \in \pi\), we have
\begin{equation*}
m_{p^{-1},p} \circ (S_p^{\mathrm{cop}} \circ S_{p^{-1}} \otimes S_{p^{-1}}) \circ \Delta_{p^{-1}} = \mu \circ \varepsilon_p \circ S_{p^{-1}} = \mu \circ \varepsilon_{p^{-1}},
\end{equation*}
\(S_p^{cop}\circ S_{p^{-1}}\) is a left inverse of \(S_{p^{-1}}\) in \((\{\operatorname{Hom}_{\mathbb{K}}(H_p, H_r)\}_{r \in \pi}, \{m_{r,s}^{p,\#}\}_{r,s \in \pi}, \mu^{p,\#})\), hence \(S_p^{cop}\circ S_{p^{-1}} = id_{p^{-1}}\)
\end{proof}

\begin{Corollary}
If \(\mathcal{H}_\pi\) is a cocommutative Hopf \(\pi\)-algebra, then for any \(p \in \pi\), \(S_p^{-1} = S_{p^{-1}}\).
\end{Corollary}
	
\section{Rota-Baxter Hopf \(\pi\)-algebras}
Based on the convolution \(\pi\)-algebra tools and antipode properties established in the previous section, this section extends Rota--Baxter operators to cocommutative Hopf \(\pi\)-algebras, presents the axiomatic definitions and weight transition relations, constructs explicit examples, and discusses the basic properties of the operators as well as the derived structures.

\begin{Definition} 
A Rota--Baxter Hopf \(\pi\)-algebra \({}^{R}\mathcal{H}_\pi\) consists of a cocommutative Hopf \(\pi\)-algebra \(\mathcal{H}_\pi\) endowed with a family of linear maps, falling into one of the following two weight classes:
\begin{enumerate}[label=(\arabic*)]
\item \textbf{Weight 1} \\
The family of maps takes the form \(\{B_{p}: H_{p} \to H_{p^{-1}}\}_{p \in \pi}\), and satisfies:
\begin{enumerate}[label=(\roman*)]
\item For all \(p \in \pi\) and \(h \in H_{p}\),
\[
B_{p}(h)_{\left( 1,p^{-1}\right) }\otimes B_{p}(h)_{\left( 2,p^{-1}\right) }=B_{p}(h_{(1,p)})\otimes B_{p}(h_{(2,p)}), \quad \varepsilon _{p^{-1}}(B_{p}(h))=\varepsilon _{p}(h).
\]
\item For any \(p, q \in \pi\), \(h \in H_{p}\) and \(g \in H_{q}\),
\[
B_{p}(h) B_{q}(g)=B_{q p}\left(h_{(1, p)} B_{p}\left(h_{(2, p)}\right) g S_{p^{-1}} B_{p}\left(h_{(3, p)}\right)\right) .
\]
\end{enumerate}

\item \textbf{Weight \(-1\)} \\
The family of maps takes the form \(\{B_{p}: H_{p} \to H_{p}\}_{p \in \pi}\), and satisfies:
\begin{enumerate}[label=(\roman*)]
\item For all \(p \in \pi\) and \(h \in H_{p}\),
\[
B_{p}(h)_{(1,p)}\otimes B_{p}(h)_{(2,p)}=B_{p}(h_{(1,p)})\otimes B_{p}(h_{(2,p)}), \quad \varepsilon _{p}(B_{p}(h))=\varepsilon _{p}(h);
\]
\item For any \(p, q \in \pi\), \(h \in H_{p}\) and \(g \in H_{q}\),
\[
B_{p}(h)B_{q}(g)=B_{pq}\Big (B_{p}(h_{(1,p)})gS_{p}B_{p}(h_{(2,p)})h_{(3,p)}\Big ).
\]
\end{enumerate}
\end{enumerate}
\end{Definition}

The two weight classes of Rota--Baxter structures are not independent of each other; a bijective correspondence between them can be established by means of the antipode operators.

\begin{Proposition}
Let \({}^{R}\mathcal{H}_\pi\) be a Rota--Baxter Hopf \(\pi\)-algebra.
\begin{enumerate}[label=(\arabic*)]
\item If \({}^{R}\mathcal{H}_\pi\) is of weight \(1\), then \((\mathcal{H}_\pi, \{R_{p}=B_{p^{-1}} \circ S_{p}\}_{p \in \pi})\) is a Rota--Baxter Hopf \(\pi\)-algebra of weight \(-1\);
\item If \({}^{R}\mathcal{H}_\pi\) is of weight \(-1\), then \((\mathcal{H}_\pi, \{B_{p}=R_{p^{-1}} \circ S_{p}\}_{p \in \pi})\) is a Rota--Baxter Hopf \(\pi\)-algebra of weight \(1\).
\end{enumerate}
\end{Proposition}

\begin{proof}
For any \(p, q \in \pi\), \(h \in H_p\) and \(g \in H_q\),
\[
\begin{aligned}
R_p(h) R_q(g) &= B_{p^{-1}} S_p(h) B_{q^{-1}} S_q(g) \\
&= B_{q^{-1}p^{-1}}\bigl(S_p(h)_{(1,p^{-1})} B_{p^{-1}}\bigl(S_p(h)_{(2,p^{-1})}\bigr) S_q(g) S_{p^{-1}} B_{p^{-1}}\bigl(S_p(h)_{(3,p^{-1})}\bigr)\bigr) \\
&= B_{q^{-1}p^{-1}}\bigl(S_p(h_{(1,p)}) B_{p^{-1}} S_p\bigl(h_{(2,p)}\bigr) S_q(g) S_p B_{p^{-1}} S_p\bigl(h_{(3,p)}\bigr)\bigr) \\
&= B_{q^{-1}p^{-1}}\bigl(S_p(h_{(1,p)}) S_{p^{-1}} S_p B_{p^{-1}}\bigl(h_{(2,p)}\bigr) S_q(g) S_p B_{p^{-1}} S_p\bigl(h_{(3,p)}\bigr)\bigr) \\
&= B_{q^{-1}p^{-1}} S_{pq}\Bigl(B_{p^{-1}} S_p\bigl(h_{(1,p)}\bigr) g \, S_p B_{p^{-1}} S_p\bigl(h_{(2,p)}\bigr) h_{(3,p)}\Bigr) \\
&= R_{pq}\Bigl(R_p\bigl(h_{(1,p)}\bigr) g \, S_p R_p\bigl(h_{(2,p)}\bigr) h_{(3,p)}\Bigr),
\end{aligned}
\]
that is, the family \(\{R_p\}_{p\in\pi}\) on the cocommutative Hopf \(\pi\)-algebra
\(\mathcal{H}_{\pi}\)
forms a Rota-Baxter algebra of weight \(-1\).

Part (2) follows by completely symmetric arguments: starting from the weight \(-1\) Rota--Baxter identity and applying the same antipode properties, one verifies that \(B_{p}=R_{p^{-1}} \circ S_{p}\) satisfies the weight 1 axioms. 
\end{proof}

\begin{Remark}
By virtue of the bijectivity of the antipode and the identity \(S_p^{-1} = S_{p^{-1}}\) from Corollary 2, the weight conversion established in Proposition 1 is a bijective correspondence between families of weight 1 and weight \(-1\) Rota--Baxter operators, and the two constructions are mutually inverse. Accordingly, all structural properties of weight 1 Rota--Baxter operators carry over fully to weight \(-1\) operators, and vice versa.
\end{Remark}

\begin{Example}
\begin{itemize}
\item[(1)]	Let \(\prescript{R}{}{\mathcal{H}}_{\pi}\) be a Rota-Baxter Hopf \(\pi\)-algebra of weight 1. Then the map \[
B : H_e \to H_e, \quad h \to S_e(h_{(1,e)}) B_e S_e(h_{(2,e)})
\] is a Rota-Baxter operator of weight 1 on the cocommutative Hopf algebra \((H_e, m_{e,e}, \mu, \Delta_e, \varepsilon_e, S_e)\).

\item[(2)] Let \(\mathcal{H}_{\pi}\) be a Rota-Baxter Hopf \(\pi\)-algebras of weight -1, then the map \[
B: H_e \to H_e,h \mapsto h_{(1,e)} B_e S_e \bigl(h_{(2,e)}\bigr)
\] is a Rota-Baxter operaters of weight \(-1\) on the cocommutative Hopf algebra \(\bigl(H_e, m_{e,e}, \mu, \Delta_e, \varepsilon_e, S_e\bigr)\).
\end{itemize}
\end{Example} 

\indent Applying Corollary 2, we have:
\begin{Proposition} 
For any cocommutative Hopf \(\pi\)-algebra \(\mathcal{H}_\pi\) we have for any \(p,q \in \pi\), \(h \in H_{pq}\),
\begin{align*}
S_{qp}(h_{(1,p)} S_p(h_{(2,p)})g S_{p^{-1}} S_p(h_{(3,p)}))
&= S_{qp}(\varepsilon_p(h_{(1,p)})g S_{p^{-1}} S_p(h_{(2,p)}))\\ 
&= S_{qp}(gS_{p}^{-1} S_p(h))
= S_{qp}(gh) = S_p(h) S_q(g),\\
\end{align*}
that is, \((\mathcal{H}_{\pi}, \{S_p\}_{p \in \pi})\) is a Rota-Baxter Hopf \(\pi\)-algebra of weight 1.
\end{Proposition} 

By virtue of the canonical Rota--Baxter property of the antipode, the most basic group-algebra type example can be constructed directly.

\begin{Example}
Let \(\pi\) be an arbitrary group. We construct a cocommutative Hopf \(\pi\)-algebra \(\mathcal{H}_\pi\) as follows.

For each \(p\in\pi\), the graded component is \(H_p = \mathbb{K}\cdot p\), the one-dimensional vector space spanned by the group element \(p\).
The graded multiplication \(m_{p,q}\colon H_p\otimes H_q \to H_{pq}\) is the linear extension of the group multiplication, satisfying
\(
m_{p,q}(p\otimes q) = pq.
\)
The unit map \(\mu\colon \mathbb{K}\to H_e\) is given by \(\mu(1_\mathbb{K}) = e\), where \(e\) denotes the identity element of \(\pi\).

For each \(p\in\pi\), the graded comultiplication \(\Delta_p\colon H_p\to H_p\otimes H_p\) and the graded counit \(\varepsilon_p\colon H_p\to \mathbb{K}\) are defined by
\(
\Delta_p(p) = p\otimes p, \varepsilon_p(p) = 1,
\)
respectively. In other words, every group element is group-like.
The antipode family \(\{S_p\colon H_p\to H_{p^{-1}}\}\) is defined by
\(
S_p(p) = p^{-1}.
\)

By proposition 3, taking \(B_p = S_p\), i.e., \(B_p(p) = p^{-1}\) for all \(p\in\pi\), yields a Rota--Baxter Hopf \(\pi\)-algebra  of weight \(1\).

Second, by proposition 2, the weight \(-1\) counterpart of the above operator is the identity family \(B_p = \mathrm{id}_{H_p}\), given by \(B_p(p) = p\) for all \(p\in\pi\).
\end{Example}

\begin{Example}
Let \(\pi = S_3\), the symmetric group on three elements. Let \(V_4 = \{e, a, b, c\}\) be the Klein four-group, with multiplication given by \(a^2 = b^2 = c^2 = e\), \(ab = ba = c\), \(bc = cb = a\), \(ca = ac = b\).

The group \(S_3\) acts on \(V_4\) by group automorphisms, permuting the three order-\(2\) elements \(a, b, c\). Identifying \(a, b, c\) with labels \(1, 2, 3\) respectively, for any \(p\in S_3\) the action satisfies
\[
p\rightharpoonup e = e,\quad p\rightharpoonup a = a_{p(1)},\quad p\rightharpoonup b = a_{p(2)},\quad p\rightharpoonup c = a_{p(3)},
\]
where \(a_1=a\), \(a_2=b\), \(a_3=c\).

We construct a crossed product cocommutative Hopf \(\pi\)-algebra \(\mathcal{H}_\pi\) as follows:
For each \(p\in\pi\), the graded component is \(H_p = \mathbb{K}V_4\).

The graded multiplication \(m_{p,q}\colon H_p\otimes H_q \to H_{pq}\) is given by the composition of group algebra multiplication and the group action, satisfying
\(m_{p,q}(h\otimes g) = h\cdot (p\rightharpoonup g)\) for any \(h,g\in V_4,\)
where \(\cdot\) denotes the multiplication in the group algebra \(\mathbb{K}V_4\).The unit map \(\mu\colon \mathbb{K}\to H_e\) satisfies \(\mu(1_\mathbb{K}) = e\).

For each \(p\in\pi\), the graded comultiplication \(\Delta_p\colon H_p\to H_p\otimes H_p\) and the graded counit \(\varepsilon_p\colon H_p\to \mathbb{K}\) satisfy
\(
\Delta_p(g) = g\otimes g,\varepsilon_p(g) = 1\) for any \(g\in V_4.\)
The antipode family \(\{S_p\colon H_p\to H_{p^{-1}}\}\) satisfies
\(
S_p(h) = p^{-1}\rightharpoonup h\)  for any \(h\in V_4.\)

Based on this Hopf \(\pi\)-algebra, taking the operator family \(B_p = S_p\), i.e., \(B_p(h) = p^{-1}\rightharpoonup h\) for all \(p\in\pi\) and \(h\in V_4\), yields a Rota--Baxter Hopf \(\pi\)-algebra  of weight \(1\), the corresponding weight \(-1\) operator is the identity family \(B_p = \mathrm{id}_{H_p}\), i.e., \(B_p(h) = h\) for all \(p\in\pi\) and \(h\in V_4\).
\end{Example}

This example focuses on the non-trivial graded multiplication induced by a non-abelian index group, with the canonical antipode Rota-Baxter operator.In contrast, we take an abelian index group below, where the crossed product Hopf \(\pi\)-algebra admits multiple distinct Rota-Baxter operators beyond the canonical one.

\begin{Example}
Let \(\pi = \mathbb{Z}/2\mathbb{Z} = \{0, 1\}\) denote the additive group of integers modulo \(2\). Let \(V_4 = \{e, a, b, c\}\) be the Klein four-group in previous example.

The group \(\pi\) acts on \(V_4\) by group automorphisms: the generator \(1\in\pi\) swaps \(a\) and \(b\) while fixing \(e\) and \(c\), i.e.,
\[
1\rightharpoonup e = e,\quad 1\rightharpoonup a = b,\quad 1\rightharpoonup b = a,\quad 1\rightharpoonup c = c.
\]

We construct a crossed-product cocommutative Hopf \(\pi\)-algebra \(\mathcal{H}_\pi\) as follows: For each \(p\in\pi\), the graded component is \(H_p = \mathbb{K}V_4\).

The graded multiplication \(m_{p,q}\colon H_p\otimes H_q \to H_{p+q}\) is defined by
\(
m_{p,q}(h\otimes g) = h\cdot (p\rightharpoonup g)\) for any \(h,g\in V_4,\)
where \(\cdot\) denotes the multiplication in the group algebra \(\mathbb{K}V_4\). The unit map \(\mu\colon \mathbb{K}\to H_0\) is given by \(\mu(1_\mathbb{K}) = e\).

For each \(p\in\pi\), the graded comultiplication \(\Delta_p\colon H_p\to H_p\otimes H_p\) and the graded counit \(\varepsilon_p\colon H_p\to \mathbb{K}\) satisfy
\(\Delta_p(h) = h\otimes h, \varepsilon_p(h) = 1\) for any  \(h\in V_4.\)
The antipode family \(\{S_p\colon H_p\to H_{-p}=H_p\}\) satisfies
\(S_p(h) = p\rightharpoonup h\) for any \(h\in V_4.\)

On this Hopf \(\pi\)-algebra, we construct three distinct families of Rota–Baxter operators. For each family, both the weight-\(1\) form and the corresponding weight-\((-1)\) form are presented below.
\begin{itemize}
\item[(1)]  The operator family of weight \(1\) is given by
\(
B_p^{(1)}(h) = p\rightharpoonup h\) for any \(h\in V_4\); the associated operator family of weight \(-1\) reads
\(
R_p^{(1)}(h) = h\) for any \(h\in V_4.\)

\item[(2)]  The operator family of weight \(1\) satisfies
\[
B_p^{(2)}(e) = e,\quad B_p^{(2)}(a) = p\rightharpoonup b,\quad B_p^{(2)}(b) = p\rightharpoonup a,\quad B_p^{(2)}(c) = c.
\]
The corresponding operator family of weight \(-1\) satisfies
\[
R_p^{(2)}(e) = e,\quad R_p^{(2)}(a) = b,\quad R_p^{(2)}(b) = a,\quad R_p^{(2)}(c) = c.
\]

\item[(3)]  The operator family of weight \(1\) satisfies
\[
B_p^{(3)}(e) = e,\quad B_p^{(3)}(a) = p\rightharpoonup c,\quad B_p^{(3)}(b) = p\rightharpoonup c,\quad B_p^{(3)}(c) = e.
\]
The corresponding operator family of weight \(-1\) satisfies
\[
R_p^{(3)}(e) = e,\quad R_p^{(3)}(a) = c,\quad R_p^{(3)}(b) = c,\quad R_p^{(3)}(c) = e.
\]
\end{itemize}
\end{Example}

The above examples verify the existence and richness of Rota--Baxter Hopf \(\pi\)-algebra structures. In what follows, we discuss the basic properties of Rota--Baxter operators from three aspects: preservation of canonical structural elements, equivalent characterization in the commutative and cocommutative case, and invariance under automorphisms of the underlying Hopf \(\pi\)-algebra.
The above examples confirm the existence and diversity of such structures. We next discuss the preservation properties of the operators with respect to typical structural elements in Hopf algebras.

\begin{Proposition}
Suppose \(\prescript{R}{}{\mathcal{H}}_{\pi}\) is a Rota-Baxter Hopf \(\pi\)-algebra of weight 1. Then: 
\begin{itemize}
\item[(1)] If \(g \in H_p\) is a group-like element in \(\mathcal{H}_p\), then \(B_p(g)\) is a group-like element in \(\mathcal{H}_{p^{-1}}\); 
\item[(2)] \(B_e \circ\mu = \mu\); 
\item[(3)] if \(a\) is a primitive element (i.e., \(a \in H_e\) and \(\Delta_e(a) = a \otimes 1 + 1 \otimes a\)), then \(B_e(a)\) is also a primitive element.
\end{itemize}
\end{Proposition}

\begin{proof}
(1) For any group-like element \(g\in H_p\), we have 
\begin{align*}\Delta_{p^{-1}} B_p(g) &= B_p(g)_{(1,p^{-1})} \otimes B_p(g)_{(2,p^{-1})}\\
& = B_p(g_{(1,p)}) \otimes B_p(g_{(2,p)}) = B_p(g) \otimes B_p(g),
\end{align*}
and
\(\varepsilon_{p^{-1}} B_p(g) = \varepsilon_p(g).\)

(2) Since \(\mu(1_K)=1\in H_e\) is a group-like element in \((H_e,\Delta_e,\varepsilon_e)\) by the conclusion of (1), \(B_e(1)\) is also a group-like element in \((H_e,\Delta_e,\varepsilon_e)\).

By the antipode compatibility condition of Hopf \(\pi\)-algebras, we have \[
B_e(1)S_e B_e(1) = \varepsilon_e B_e(1)1 = 1= \varepsilon_e B_e(1)1=S_e B_e(1)B_e(1).
\] 

From the defining identity of a Rota-Baxter Hopf \(\pi\)-algebra, we have \[
B_e(1)B_e(1) = B_e(1B_e(1)1S_e B_e(1))1 = B_e(1).
\] 

Multiplying both sides on the left by \(S_e B_e(1)\), we obtain \[
B_e(1)=S_e B_e(1)B_e(1)B_e(1)=S_e B_e(1) B_e(1) =1.
\] 

(3) For any primitive element \(a\in H_e\), we have
\begin{align*}
\Delta_e B_e(a) &= B_e(a)_{(1,e)} \otimes B_e(a)_{(2,e)}\\
& = B_e(a_{(1,e)}) \otimes B_e(a_{(2,e)})= (B_e \otimes B_e)\big[(a \otimes 1)+(1 \otimes a)\big]\\
&= B_e(a) \otimes B_e(1) + B_e(1) \otimes B_e(a)= B_e(a) \otimes 1 + 1 \otimes B_e(a).
\end{align*} 
\end{proof}
	
\begin{Proposition}
If \(\mathcal{H}_\pi\) is a cocommutative and commutative Hopf \(\pi\)-algebra, and \(\{B_p : H_p \to H_{p^{-1}}\}_{p \in \pi}\) is a family of coalgebra homomorphisms. Then \(\prescript{R}{}{\mathcal{H}}_{\pi}\) is a Rota-Baxter Hopf \(\pi\)-algebra of weight 1 if and only \(B_{pq}\circ m_{p,q} = m_{p^{-1},q^{-1}} \circ(B_p \otimes B_q)\) for any \(p, q \in \pi\).
\end{Proposition}
	
\begin{proof}
For any \(p, q \in \pi\), \(h \in H_p\), \(g \in H_q\), we have
\begin{align*}
B_{qp}\big(h_{(1,p)}B_p(h_{(2,p)})gS_{p^{-1}}B_p(h_{(3,p)})\big)&= B_{qp}\big(h_{(1,p)}B_p(h_{(2,p)})S_{p^{-1}}B_p(h_{(3,p)})g\big) \\
&= B_{qp}\big(h_{(1,p)}B_p(h_{(21,p)})S_{p^{-1}}B_p(h_{(22,p)})g\big)\\
&= B_{qp}\big(h_{(1,p)}B_p(h_{(2,p)})_{(1,p^{-1})}S_{p^{-1}}\big[B_p(h_{(2,p)})_{(2,p^{-1})}\big]g\big) \\
&= B_{qp}\big(h_{(1,p)}\varepsilon_{p^{-1}}B_p(h_{(2,p)})g\big)\\
&= B_{qp}\big(h_{(1,p)}\varepsilon_p(h_{(2,p)})g\big)
= B_{qp}(hg)
= B_{pq}(hg).
\end{align*}
\end{proof}

\indent Applying Theorem 5, we have:
\begin{Proposition}
Suppose \(\prescript{R}{}{\mathcal{H}}_{\pi}\) is a Rota-Baxter Hopf \(\pi\)-algebra of weight 1 and \(\{\varphi_p\}_{p \in \pi}\) is a Hopf \(\pi\)-algebra automorphism or anti-automorphism of \(\mathcal{H}_\pi\). Let \(\{B'_p = \varphi_{p^{-1}} \circ B_p \circ \varphi_p^{-1}\}_{p \in \pi}\). Then \((\mathcal{H}_{\pi}, \{B'_p\}_{p \in \pi})\) is also a Rota-Baxter Hopf \(\pi\)-algebra of weight 1.
\end{Proposition}
	
\begin{proof}
When \(\{\varphi_p\}_{p \in \pi}\) is an automorphism, we have
\begin{align*}
B_p'(h)B_q'(g)
&= \varphi_{p^{-1}} B_p \varphi_p^{-1}(h)\varphi_{q^{-1}} B_q \varphi_q^{-1}(g)
= \varphi_{p^{-1}q^{-1}} \big(B_p \varphi_p^{-1}(h)B_q \varphi_q^{-1}(g)\big) \\
&= \varphi_{p^{-1}q^{-1}} B_{qp} \big(\varphi_p^{-1}(h)_{(1,p)} B_p \big(\varphi_p^{-1}(h)_{(2,p)}\big)\varphi_q^{-1}(g)S_{p^{-1}} B_p \big(\varphi_p^{-1}(h)_{(3,p)}\big)\big) \\
&= \varphi_{p^{-1}q^{-1}} B_{qp} \big(\varphi_p^{-1}(h_{(1,p)}) B_p \big(\varphi_p^{-1}(h_{(2,p)})\big)\varphi_q^{-1}(g)S_{p^{-1}} B_p \big(\varphi_p^{-1}(h_{(3,p)})\big)\big) \\
&= \varphi_{q^{-1}p^{-1}} B_{qp} \big(\varphi_p^{-1}(h_{(1,p)}) B_p \varphi_p^{-1}(h_{(2,p)})\varphi_q^{-1}(g)S_{p^{-1}} B_p \varphi_p^{-1}(h_{(3,p)})\big) \\
&= \varphi_{q^{-1}p^{-1}} B_{qp} \big(\varphi_p^{-1}(h_{(1,p)}) B_p \varphi_p^{-1}(h_{(2,p)})\varphi_q^{-1}(g)\varphi_p^{-1}\varphi_p S_{p^{-1}} B_p \varphi_p^{-1}(h_{(3,p)})\big) \\
&= \varphi_{q^{-1}p^{-1}} B_{qp} \big(\varphi_p^{-1}(h_{(1,p)}) B_p \varphi_p^{-1}(h_{(2,p)})\varphi_q^{-1}(g)\varphi_p^{-1}S_{p^{-1}}\varphi_{p^{-1}} B_p \varphi_p^{-1}(h_{(3,p)})\big) \\
&= \varphi_{p^{-1}q^{-1}} B_{qp} \varphi_{qp}^{-1}\big(h_{(1,p)}\big)\varphi_{p^{-1}} B_p \varphi_p^{-1}(h_{(2,p)})gS_{p^{-1}}\varphi_{p^{-1}} B_p \varphi_p^{-1}(h_{(3,p)})\big)\\
&=B_{qp}'\big(h_{(1,p)}\big)B_p'\big(h_{(2,p)}\big)gS_{p^{-1}}B_p'\big(h_{(3,p)}\big). 
\end{align*}
The case of anti-automorphism can be proved similarly.
\end{proof}

Furthermore, Rota--Baxter operators can also induce new graded multiplications and antipode families, giving rise to a class of derived Hopf \(\pi\)-algebras.

\begin{Theorem}
Let \({}^{R}\mathcal{H}_\pi\) be a Rota‑Baxter Hopf \(\pi\)-algebra. Define two families of linear maps \[
\bigl\{m_{p,q}^{*}\colon H_{p}\otimes H_{q}\to H_{qp},\allowbreak\ m_{p,q}^{*}(h\otimes g)= h *_{p,q} g =h_{(1,p)}B_{p}(h_{(2,p)})gS_{p^{-1}}\bigl(B_{p}(h_{(3,p)})\bigr)\bigr\}_{p,q\in\pi}
\] and \[
\bigl\{T_{p}\colon H_{p}\to H_{p^{-1}},\allowbreak\ T_{p}(h)=S_{p^{-1}}\bigl(B_{p}(h_{(1,p)})\bigr)S_{p}(h_{(2,p)})B_{p}(h_{(3,p)})\bigr\}_{p\in\pi}
\].
Then the following statements hold:
\begin{enumerate}
\item[(1)]For any \(p, q \in \pi\), \(h \in H_p\) and \(g \in H_q\), the following identities hold:
\[
\Delta_{qp}\left(h *_{p,q} g\right) = \left(h_{(1,p)} *_{p,q} g_{(1,q)}\right) \otimes \left(h_{(2,p)} *_{p,q} g_{(2,q)}\right),\quad
\varepsilon_{qp}\left(h *_{p,q} g\right) = \varepsilon_p(h)\varepsilon_q(g).
\]

\item[(2)] \allowbreak For any \(p \in \pi\) and \(h \in H_p\), the following identity holds:
\[
\varepsilon_p(h) 1 = B_p(h_{(1,p)}) \, B_{p^{-1}}\big(T_p(h_{(2,p)})\big).
\]

\item[(3)] \(\left(\{H_p\}_{p\in\pi}, \{m_{p,q}^*\}_{p,q\in\pi}, \mu\right)\) is a \(\pi^{\mathrm{op}}\)-algebra.

\item[(4)] \(\bigl(\{H_p, \Delta_p, \varepsilon_p\}_{p\in\pi},\allowbreak \{m^*_{p,q}\}_{p,q\in\pi},\allowbreak \mu,\allowbreak \{T_p\}_{p\in\pi}\bigr)\) is a cocommutative Hopf \(\pi^{op}\)-algebra, called the \emph{descendant Hopf \(\pi\)-algebra} of \(\mathcal{H}_{\pi}\) with respect to \(\{B_p\}_{p\in\pi}\).

\item[(5)] \(\{B_p\}_{p\in\pi}\) is a generalized graded Hopf algebra homomorphism from \(\bigl(\{H_p, \Delta_p, \varepsilon_p\}_{p\in\pi},\allowbreak \{m^*_{p,q}\}_{p,q\in\pi},\allowbreak \mu,\allowbreak \{T_p\}_{p\in\pi}\bigr)\) to \(\mathcal{H}_{\pi}\) .

\item[(6)]   \(\bigl(\{H_p, \Delta_p, \varepsilon_p\}_{p\in\pi},\allowbreak \{m^*_{p,q}\}_{p,q\in\pi},\allowbreak \mu,\allowbreak \{T_p\}_{p\in\pi},\{B_p\}_{p\in\pi}\bigr)\) is a Rota‑Baxter Hopf \(\pi\)-algebra.
\end{enumerate}
	
\begin{proof}
(1)For any \(p,q\in\pi\), \(h\in H_p\), \(g\in H_q\),we have
\[
\begin{aligned}
&\Delta_{qp}\left(h *_{p,q} g\right)= \Delta_{qp}\left( h_{(1,p)} B_p(h_{(2,p)}) g\, S_{p^{-1}}\big(B_p(h_{(3,p)})\big) \right) \\
&= \left( h_{(11,p)} B_p(h_{(21,p)}) g_{(1,q)} S_{p^{-1}}\big(B_p(h_{(31,p)})\big) \right)
\otimes \left( h_{(12,p)} B_p(h_{(22,p)}) g_{(2,q)} S_{p^{-1}}\big(B_p(h_{(32,p)})\big) \right) \\
&= \left( h_{(11,p)} B_p(h_{(12,p)}) g_{(1,q)} S_{p^{-1}}\big(B_p(h_{(13,p)})\big) \right) 
\otimes \left( h_{(21,p)} B_p(h_{(22,p)}) g_{(2,q)} S_{p^{-1}}\big(B_p(h_{(23,p)})\big) \right) \\
&= \left( h_{(1,p)} *_{p,q} g_{(1,q)} \right) \otimes \left( h_{(2,p)} *_{p,q} g_{(2,q)} \right)
\end{aligned}
\]
and
\[
\begin{aligned}
\varepsilon_{qp}\left(h *_{p,q} g\right)
&= \varepsilon_{qp}\left( h_{(1,p)} B_p(h_{(2,p)}) g\, S_{p^{-1}}\big(B_p(h_{(3,p)})\big) \right) \\
&= \varepsilon_p(h_{(1,p)})\, \varepsilon_{p^{-1}}\!\big(B_p(h_{(2,p)})\big)\, \varepsilon_q(g)\, \varepsilon_{p}S_{p^{-1}}(B_p(h_{(3,p)})) \\
&= \varepsilon_p(h_{(1,p)})\, \varepsilon_p(h_{(2,p)})\, \varepsilon_q(g)\, \varepsilon_p(h_{(3,p)})= \varepsilon_p(h)\,\varepsilon_q(g).
\end{aligned}
\]

(2)For any \(p\in\pi\) and \(h\in H_p\), we have
\[
\begin{aligned}
&B_p(h_{(1,p)}) \, B_{p^{-1}}\big(T_p(h_{(2,p)})\big)= B_e(h_{(11,p)}\, B_p(h_{(12,p)}) \, T_p(h_{(2,p)}) \, S_{p^{-1}}\big(B_p(h_{(3,p)}))\big) \\
&= B_e\Big( h_{(11,p)} B_p(h_{(12,p)}) S_{p^{-1}}\big(B_p(h_{(21,p)})\big) S_{p}\big(h_{(22,p)})\big)B_p(h_{(23,p)}) S_{p^{-1}}\big(B_p(h_{(3,p)})\big) \Big)\\
&= B_e\Big( h_{(1,p)} B_p(h_{(21,p)}) S_{p^{-1}}\big(B_p(h_{(22,p)})\big) S_{p}\big(h_{(3,p)})\big)B_p(h_{(41,p)}) S_{p^{-1}}\big(B_p(h_{(42,p)})\big) \Big)\\
&= B_e\Big( h_{(1,p)} B_p(h_{(2,p)})_{(1,p^{-1})} S_{p^{-1}}\big(B_p(h_{(2,p)})_{(2,p^{-1})}\big) S_{p}\big(h_{(3,p)})\big)B_p(h_{(4,p)})_{(1,p^{-1})}\\ \allowbreak &S_{p^{-1}}\big(B_p(h_{(4,p)})_{(2,p^{-1})}\big) \Big)\\
&= B_e\big(h_{(1,p)} \, \varepsilon_{p^{-1}}B_p(h_{(2,p)}) \, S_p(h_{(3,p)}) \, \varepsilon_{p^{-1}}B_p(h_{(4,p)})\big) \\
&= B_e\big(h_{(1,p)}) \, \varepsilon_p(h_{(2,p)}) \, S_p(h_{(3,p)}) \, \varepsilon_p(h_{(4,p)}) \big)\\
&= B_e\big(h_{(1,p)}) \, S_p(h_{(2,p)})\big)
= \varepsilon_p(h) \, B_e(1)
= \varepsilon_p(h)1.
\end{aligned}
\]

(3) For any \(p,q\in\pi\), \(h\in H_p\), \(g\in H_q\), by the original Rota–Baxter identity,
\[
B_{pq}\bigl(m^*_{p,q}(h\otimes g)\bigr) = B_p(h)B_q(g),
\]
thus \(\{B_p\}_{p\in\pi}\) preserves multiplications.For any \(p,q,r\in\pi\), \(h\in H_p\), \(g\in H_q\), \(j\in H_r\), we have
\[
\begin{aligned}
&\big( m_{qp, r}^* \circ ( m_{p,q}^* \otimes \mathrm{id}_{H_r} ) \big) (h \otimes g \otimes j)
= \big( h *_{p,q} g \big) *_{qp, r} j \\
&= \big(h *_{p,q} g\big)_{(1,qp)} B_{qp}\!\left( \big(h *_{p,q} g\big)_{(2,qp)} \right) j \, S_{(qp)^{-1}} B_{qp}\!\left( \big(h *_{p,q} g\big)_{(3,qp)} \right)\\
&= \big(h_{(1,p)} *_{p,q} g_{(1,q)}\big) \, B_{qp}\!\big(h_{(2,p)} *_{p,q} g_{(2,q)}\big) \, j \, S_{(qp)^{-1}} B_{qp}\!\big(h_{(3,p)} *_{p,q} g_{(3,q)}\big)\\
&= h_{(11,p)} B_p\big(h_{(12,p)}\big) g_{(1,q)} S_{p^{-1}} B_p\big(h_{(13,p)}\big) B_p\big(h_{(2,p)}\big) B_q\big(g_{(2,q)}\big) j S_{(qp)^{-1}} \big( B_p\big(h_{(3,p)}\big) B_q\big(g_{(3,q)}\big) \big)\\
&= h_{(1,p)} B_p\big(h_{(2,p)}\big) g_{(1,q)} S_{p^{-1}} B_p\big(h_{(3,p)}\big) B_p\big(h_{(4,p)}\big) B_q\big(g_{(2,q)}\big) j S_{q^{-1}}B_q(g_{(3,q)}) S_{p^{-1}}B_p(h_{(5,p)})\\
&= h_{(1,p)} B_p\big(h_{(2,p)}\big) g_{(1,q)} \varepsilon(h_{(3,p)})  B_q\big(g_{(2,q)}\big) j S_{q^{-1}}B_q(g_{(3,q)}) S_{p^{-1}}B_p(h_{(4,p)})\\
&= h_{(1,p)} B_p\big(h_{(2,p)}\big) g_{(1,q)}   B_q\big(g_{(2,q)}\big) j S_{q^{-1}}B_q(g_{(3,q)}) S_{p^{-1}}B_p(h_{(3,p)})\\
&= h_{(1,p)} B_p\big(h_{(2,p)}\big) \big( g *_{q,r} j \big) S_{p^{-1}} B_p\big(h_{(3,p)}\big)\\
&= h *_{p, rq} \big( g *_{q,r} j \big) =\big( m_{p, rq}^* \circ \big( \mathrm{id}_{H_p} \otimes m_{q,r}^* \big) \big) (h \otimes g \otimes j)
\end{aligned}
\]
and
\[
\begin{aligned}
m_{e,p}^*(h\otimes 1)&=h *_{p,e} 1 = h_{(1,p)} B_p\big(h_{(2,p)}\big) S_{p^{-1}} B_p\big(h_{(3,p)}\big)= h_{(1,p)} \varepsilon_{p^{-1}} B_p\big(h_{(2,p)}\big)\\
& = h_{(1,p)} \varepsilon_p\big(h_{(2,p)}\big)=h =1 B_e(1) h S_e B_e(1)=1 *_{e,p} h= m_{p,e}^*(1\otimes h). 
\end{aligned}
\]

(4)	For any \(p\in\pi\), \(h\in H_p\),
\[		
\begin{aligned}
(m_{p,p^{-1}}^* &\circ (\mathrm{id}_{H_p}\otimes T_p ) \circ \Delta_p )(h)=h_{(1,p)} *_{p,p^{-1}} T_p(h_{(2,p)})\\
&=h_{(11,p)} B_p\big(h_{(12,p)}\big) T_p\big(h_{(2,p)}\big) S_{p^{-1}} B_p\big(h_{(13,p)}\big)\\
&=h_{(11,p)} B_p \big( h_{(12,p)} \big) S_{p^{-1}} B_p \big( h_{(21,p)} \big) S_p \big( h_{(22,p)} \big) B_p \big( h_{(23,p)} \big) S_{p^{-1}} B_p \big( h_{(13,p)} \big)\\
&=h_{(1,p)} B_p \big( h_{(2,p)} \big) S_{p^{-1}} B_p \big( h_{(3,p)} \big) S_p \big( h_{(4,p)} \big) B_p \big( h_{(5,p)} \big) S_{p^{-1}} B_p \big( h_{(6,p)} \big)\\
&=h_{(1,p)} \varepsilon_{p^{-1}} B_p \big( h_{(2,p)} \big) S_p \big( h_{(3,p)} \big) \varepsilon_{p^{-1}} B_p \big( h_{(4,p)} \big)\\
&= h_{(1,p)} \varepsilon_p \big( h_{(2,p)} \big) S_p \big( h_{(3,p)} \big) \varepsilon_p \big( h_{(4,p)} \big)= h_{(1,p)} S_p \big( h_{(2,p)} \big) = \varepsilon_p(h) 1= (\mu \circ \varepsilon_p)(h)
\end{aligned}
\]

By
\[
B_p\big(h_{(1,p)}\big) B_{p^{-1}}\big( T_p\big(h_{(2,p)}\big) \big)= B_e\big( h_{(1,p)} *_{p,p^{-1}} T_p\big(h_{(2,p)}\big) \big)
= B_e\big( \varepsilon_p(h) 1 \big) = \varepsilon_p(h) 1,
\]
we have
\[
\begin{aligned}
B_{p^{-1}}\big( T_p\big(h_{(1,p)}\big) \big) B_p\big(h_{(2,p)}\big)&= B_{p^{-1}} T_p \big(h_{(1,p)}\big) \varepsilon_{p^{-1}} B_p \big(h_{(2,p)}\big) B_p \big(h_{(3,p)}\big)\\
&= \varepsilon_{p^{-1}} B_p \big(h_{(1,p)}\big) B_{p^{-1}} T_p \big(h_{(2,p)}\big) B_p \big(h_{(3,p)}\big)\\
&= S_{p^{-1}} B_p \big(h_{(1,p)}\big) B_p \big(h_{(2,p)}\big) B_{p^{-1}} T_p \big(h_{(3,p)}\big) B_p \big(h_{(4,p)}\big)\\
&= S_{p^{-1}} B_p \big(h_{(1,p)}\big) \varepsilon_p \big(h_{(2,p)}\big) B_p \big(h_{(3,p)}\big)= \varepsilon_p(h) 1.
\end{aligned}
\]

\[
\begin{aligned}
( m_{p^{-1},p}^* &\circ ( T_p \otimes \mathrm{id}_{H_p} ) \circ \Delta_p ) (h)
=T_p\big(h_{(1,p)}\big) *_{p^{-1},p} h_{(2,p)} \\
&= T_p\big(h_{(1,p)}\big) B_{p^{-1}} T_p\big(h_{(2,p)}\big) h_{(3,p)} S_p B_{p^{-1}} T_p\big(h_{(4,p)}\big) \\
&= S_{p^{-1}} B_p\big(h_{(1,p)}\big) S_p\big(h_{(2,p)}\big) B_p\big(h_{(3,p)}\big) B_{p^{-1}} T_p\big(h_{(4,p)}\big) h_{(5,p)} S_p B_{p^{-1}} T_p\big(h_{(6,p)}\big)\\
&= S_{p^{-1}} B_p \big(h_{(1,p)}\big) S_p \big(h_{(2,p)}\big) \varepsilon_p \big(h_{(3,p)}\big) h_{(4,p)} S_p B_{p^{-1}} T_p \big(h_{(5,p)}\big) \\
&= S_{p^{-1}} B_p \big(h_{(1,p)}\big) S_p \big(h_{(2,p)}\big) h_{(3,p)} S_p B_{p^{-1}} T_p \big(h_{(4,p)}\big) \\
&= S_{p^{-1}} B_p \big(h_{(1,p)}\big) S_p B_{p^{-1}} T_p \big(h_{(2,p)}\big) \\
&= S_e \big( B_{p^{-1}} T_p \big(h_{(2,p)}\big) B_p \big(h_{(1,p)}\big) \big) = \varepsilon_p(h) 1=(\mu \circ \varepsilon_p)(h)
\end{aligned}
\]
Therefore \(\bigl(\{H_p, \Delta_p, \varepsilon_p\}_{p\in\pi},\allowbreak \{m^*_{p,q}\}_{p,q\in\pi},\allowbreak \mu,\allowbreak \{T_p\}_{p\in\pi}\bigr)\) is a Hopf \(\pi^{\mathrm{op}}\)-algebra.

(5) Clearly, \(\Phi \colon \pi^{\mathrm{op}} \to \pi\) given by \(p \mapsto p^{-1}\) is a group homomorphism. For any \(p,q \in \pi\), \(h \in H_p\) and \(g \in H_q\), we have
\[
\bigl(m_{\Phi(p),\Phi(q)} \circ (B_p \otimes B_q)\bigr)(h\otimes g)
= B_p(h) B_q(g)
= B_{pq}\bigl(h *_{p,q} g\bigr)
= \bigl(B_{pq} \circ m_{p,q}^*\bigr)(h\otimes g),
\]
which means that \(\{B_p\}_{p\in\pi}\) is a generalized graded algebra homomorphism from \(\bigl(\{H_p\}_{p\in\pi}, \{m_{p,q}^*\}_{p,q\in\pi}, \mu\bigr)\) to \(H_\pi\).

Combining with Corollary 1, the proposition is proved.

(6) For any \(p,q \in \pi\), \(h \in H_p\) and \(g \in H_q\),
\begin{equation*}
\begin{aligned}
&B_{qp} \big( h_{(1,p)} * B_p (h_{(2,p)}) * g * S_{p^{-1}} B_p (h_{(3,p)}) \big) \\
=\ & B_p (h_{(1,p)}) B_{p^{-1}} \big( B_p (h_{(2,p)}) \big) B_q (g) B_p S_{p^{-1}} B_p (h_{(3,p)}) \\
=\ & B_p (h_{(1,p)}) B_{p^{-1}} B_p (h_{(2,p)}) B_q (g) S_p B_{p^{-1}} B_p (h_{(3,p)}) \\
=\ & B_p (h) *_{p^{-1}, q^{-1}} B_q (g)
\end{aligned}
\end{equation*}
\end{proof}
\end{Theorem}

As a complementary result, the image set of the Rota--Baxter operators constitutes a Hopf \(\pi\)-subalgebra of the original Hopf \(\pi\)-algebra.

\begin{Proposition}
Let \({}^{R}\mathcal{H}_\pi\) be a  Rota--Baxter Hopf \(\pi\)-algebra, and set \(G_p = \operatorname{Im} B_{p^{-1}} \subseteq H_p\). Then \(\bigl(\{G_p\}_{p \in \pi}, \Delta_{p}, \varepsilon_{p}\}_{p \in \pi},\{m_{p, q}\}_{p, q \in \pi}, \mu, \{S_{p}\}_{p \in \pi}\bigr)\) forms a Hopf \(\pi\)-subalgebra of \(\mathcal{H}_{\pi}\).
\end{Proposition}

\begin{proof}
For any \(p, q \in \pi\), 
\(G_p \cdot G_q=(B_{p^{-1}}(H_{p^{-1}}) \cdot B_{q^{-1}}(H_{q^{-1}}) \subset B_{q^{-1}p^{-1}}(H_{q^{-1}p^{-1}}) = B_{(pq)^{-1}}(H_{(pq)^{-1}})=G_{pq}\). 
By Proposition 1, \(B_e(1) = 1\in G_e\). Hence \(\bigl(\{G_p\}_{p \in \pi},\{m_{p, q}\}_{p, q \in \pi}, \mu\bigr)\)is a \(\pi\)-subalgebra of \(H_{\pi}\).

By the coalgebra homomorphism property of \(B_{p^{-1}}\), 
\small \(\Delta_p(G_p) = \Delta_p\bigl(B_{p^{-1}}(H_{p^{-1}})\bigr) = (B_{p^{-1}} \otimes B_{p^{-1}}) \Delta_{p^{-1}}(H_{p^{-1}}) \subseteq G_p \otimes G_p\), \normalsize
and \(\varepsilon_p(G_p) = \varepsilon_{p^{-1}}(H_{p^{-1}}) \in \mathbb{K}\). Thus \(\{G_p\}\) is a \(\pi\)-subcoalgebra.

By Theorem 6 (5), \(B_p \circ T_{p^{-1}} = S_p \circ B_{p^{-1}}\) for any \(p\in \pi\). Therefore \(S_p(G_p) = S_p(B_{p^{-1}}(H_{p^{-1}})) = B_p(T_{p^{-1}}(H_{p^{-1}})) \subseteq \operatorname{Im} B_p = G_{p^{-1}}\).
\end{proof}

\section{Conclusion}

This paper presents a systematic generalization of the theory of Rota--Baxter operators from classical cocommutative Hopf algebras to the group-graded setting of Hopf \(\pi\)-algebras. The development proceeds in three stages: construction of graded convolution algebras, derivation of core structural properties of Hopf \(\pi\)-algebras, and establishment of the graded Rota--Baxter operator framework.

First, as the central technical tool of the paper, we introduce the notion of a convolution \(\pi\)-algebra induced by a coassociative counital coalgebra and a \(\pi\)-algebra, and verify that it satisfies the axioms of a graded algebra. We further show that coalgebra homomorphisms pull back to graded homomorphisms of convolution \(\pi\)-algebras, and that generalized graded homomorphisms of \(\pi\)-algebras push forward to generalized graded homomorphisms of the corresponding convolution \(\pi\)-algebras. These results provide a unified convolution calculus for the study of Hopf \(\pi\)-algebras.

Second, using convolution \(\pi\)-algebras and the notion of convolution inverse, we derive a series of fundamental properties of Hopf \(\pi\)-algebras. We give explicit formulas for the convolution inverses of two canonical classes of maps, and deduce the antipode compatibility criterion for generalized homomorphisms between Hopf \(\pi\)-algebras. We prove that the antipode of a Hopf \(\pi\)-algebra is both an anti-homomorphism of graded algebras and an anti-homomorphism of coalgebras. For cocommutative Hopf \(\pi\)-algebras, we establish the bijectivity of the antipode and show that the co-opposite graded structure again admits a Hopf \(\pi\)-algebra structure. These results complement and refine the foundational theory of Hopf \(\pi\)-algebras.

Third, based on the antipode and homomorphism properties established above, we develop the axiomatic theory of Rota--Baxter operators on cocommutative Hopf \(\pi\)-algebras, and introduce two families of graded Rota--Baxter operators, of weight \(1\) and weight \(-1\) respectively. We prove that the two families are in bijective correspondence via the antipode, so that all structural properties transfer between them. We then construct several explicit examples, including the trivial group-graded Hopf \(\pi\)-algebra, crossed product Hopf \(\pi\)-algebras with non-abelian index group \(S_3\), and crossed product Hopf \(\pi\)-algebras with abelian index group \(\mathbb{Z}/2\mathbb{Z}\); the last example yields three distinct families of Rota--Baxter operators, demonstrating the richness of the theory. We also investigate basic properties of the operators: preservation of canonical structural elements such as group-like elements and primitive elements, an equivalent characterization in the commutative and cocommutative case, and invariance under automorphisms and anti-automorphisms of the underlying Hopf \(\pi\)-algebra. Finally, we describe two derived structures induced by a Rota--Baxter operator: the descendant Hopf \(\pi^\mathrm{op}\)-algebra and the image Hopf \(\pi\)-subalgebra.

This work extends the classical correspondence between Rota--Baxter operators and Hopf algebras to the group-graded setting, and enriches the repertoire of operator structures on Hopf \(\pi\)-algebras. Building on the results of this paper, we are currently generalizing the framework to multiplier Hopf \(\pi\)-algebras. We expect that future research will elucidate the intrinsic relation between graded Rota--Baxter structures and graded Yang--Baxter equations, extending the scope of the classical results.

\section*{Funding}
This work was supported by the National Natural Science Foundation of China (Grant No. 12505111).

\section*{Acknowledgments}
The authors appreciate the forthcoming valuable comments and suggestions from the journal editor and anonymous referee, which will further refine the theoretical content and presentation of this article.

\section*{Author Contributions}
All authors jointly completed the theoretical derivation and example construction of Rota-Baxter Hopf \(\pi\)-algebras. G.S. proposed the research idea and guided the whole research work. X.T. and X.W. finished the main theorem proofs and drafted the original manuscript. Y.Z. sorted out the references and participated in the revision of the full text. All authors have read and agreed to the final published version of the manuscript.

\section*{Data Availability Statement}
No new data were created.

\section*{Conflict of Interest}
The authors declare no conflicts of interest.



\begin{thebibliography}{99}

\bibitem{baxter1960analytic}
Baxter, G. An analytic problem whose solution follows from a simple algebraic identity.
\textit{Pac. J. Math.} \textbf{1960}, \textit{10}, 731--742.

\bibitem{rota1969baxter}
Rota, G.C. Baxter algebras and combinatorial identities I, II.
\textit{Bull. Am. Math. Soc.} \textbf{1969}, \textit{75}, 325--334.

\bibitem{cartier1972structure}
Cartier, P. On the structure of free Baxter algebras.
\textit{Adv. Math.} \textbf{1972}, \textit{9}, 253--265.

\bibitem{connes2000renormalization}
Connes, A.; Kreimer, D. Renormalization in quantum field theory and the Riemann-Hilbert problem. I. The Hopf algebra structure of graphs and the main theorem.
\textit{Commun. Math. Phys.} \textbf{2000}, \textit{210}, 249--273.

\bibitem{andrews2003baxter}
Andrews, G.E.; Guo, L.; Keigher, W. Baxter algebras and Hopf algebras.
\textit{Trans. Am. Math. Soc.} \textbf{2003}, \textit{355}, 4639--4656.

\bibitem{guo2012introduction}
Guo, L. \textit{An Introduction to Rota-Baxter Algebra};
International Press: Somerville, MA, USA, 2012.

\bibitem{hopf1941topologie}
Hopf, H. Über die Topologie der Gruppen-Mannigfaltigkeiten und ihrer Verallgemeinerungen.
\textit{Ann. Math.} \textbf{1941}, \textit{42}, 22--52.

\bibitem{milnor1965structure}
Milnor, J.W.; Moore, J.C. On the structure of Hopf algebras.
\textit{Ann. Math.} \textbf{1965}, \textit{81}, 211--264.

\bibitem{sweedler1969hopf}
Sweedler, M.E. \textit{Hopf Algebras};
W.A. Benjamin: New York, NY, USA, 1969.

\bibitem{radford2012hopf}
Radford, D.E. \textit{Hopf Algebras};
World Scientific: Singapore, 2012.

\bibitem{turaev2000homotopy}
Turaev, V.G. Homotopy field theory in dimension 3 and group-categories.
arXiv:math/0012073, 2000.

\bibitem{turaev2010homotopy}
Turaev, V.G. \textit{Homotopy Quantum Field Theory};
EMS Tracts in Mathematics, vol. 10; European Mathematical Society: Zürich, 2010.

\bibitem{ohtsuki1993colored}
Ohtsuki, T. Colored ribbon Hopf algebras and universal invariants of framed links.
\textit{J. Knot Theory Ramif.} \textbf{1993}, \textit{2}, 211--232.

\bibitem{virelizier2002hopf}
Virelizier, A. Hopf group-coalgebras.
\textit{J. Pure Appl. Algebra} \textbf{2002}, \textit{171}, 75--122.

\bibitem{goncharov2021rota}
Goncharov, M. Rota-Baxter operators on cocommutative Hopf algebras.
\textit{J. Algebra} \textbf{2021}, \textit{582}, 39--56.

\bibitem{ma2022rotabaxter}
Ma, T.; Li, J.; Chen, L.; Wang, S. Rota-Baxter operators on Turaev’s Hopf group (co)algebras I: Basic definitions and related algebraic structures.
\textit{J. Geom. Phys.} \textbf{2022}, \textit{175}, 104469.

\end{thebibliography}
\end{document}